\documentclass[11pt,a4paper,twoside]{article}
\usepackage[T1]{fontenc}
\usepackage{lmodern,amsmath,amssymb,amsthm,mathtools}
\usepackage{booktabs,array,microtype,enumitem,fancyhdr,mathrsfs}
\usepackage[a4paper,inner=30mm,outer=30mm,top=27mm,bottom=26mm,headheight=14pt,headsep=13pt]{geometry}
\usepackage[hidelinks,pdfencoding=auto,psdextra]{hyperref}
\usepackage{bookmark,embedfile}
\newcommand{\papertitle}{Commuting pairs in lattices of reductive Lie algebras}
\hypersetup{pdftitle={Commuting pairs in lattices of reductive Lie algebras},
 pdfauthor={Caihua Luo, Xu Song, Yilin Zeng},
 pdfkeywords={Commuting pairs, Lie lattices, exponential sums, integral points}}
\setlist[enumerate]{topsep=4pt,itemsep=2pt}
\numberwithin{equation}{section}
\newtheorem{theorem}{Theorem}[section]
\newtheorem{corollary}[theorem]{Corollary}
\newtheorem{proposition}[theorem]{Proposition}
\newtheorem{lemma}[theorem]{Lemma}
\theoremstyle{remark}
\newtheorem*{remark}{Remark}
\newcommand{\Q}{\mathbb Q}
\newcommand{\Z}{\mathbb Z}
\newcommand{\R}{\mathbb R}
\newcommand{\C}{\mathbb C}
\newcommand{\F}{\mathbb F}
\newcommand{\K}{\mathbb K}
\newcommand{\A}{\mathbb A}
\newcommand{\PP}{\mathbb P}
\newcommand{\OO}{\mathcal O}
\newcommand{\g}{\mathfrak g}
\newcommand{\s}{\mathfrak s}
\newcommand{\z}{\mathfrak z}
\newcommand{\h}{\mathfrak h}
\newcommand{\tor}{\mathfrak t}
\newcommand{\p}{\mathfrak p}
\newcommand{\id}{\mathfrak a}
\newcommand{\diff}{\mathfrak d}
\newcommand{\lev}{\mathfrak l}
\newcommand{\slie}{\mathfrak{sl}}
\newcommand{\glie}{\mathfrak{gl}}
\newcommand{\gsp}{\mathfrak{gsp}}
\newcommand{\gso}{\mathfrak{gso}}
\newcommand{\spLie}{\mathfrak{sp}}
\newcommand{\soLie}{\mathfrak{so}}
\newcommand{\GL}{\mathrm{GL}}
\newcommand{\Res}{\operatorname{Res}}
\newcommand{\vol}{\operatorname{vol}}
\newcommand{\Tr}{\operatorname{Tr}}
\newcommand{\Ad}{\operatorname{Ad}}
\newcommand{\ad}{\operatorname{ad}}
\newcommand{\Lie}{\operatorname{Lie}}
\newcommand{\Av}{\operatorname{Av}}
\newcommand{\sinc}{\operatorname{sinc}}
\newcommand{\covol}{\operatorname{covol}}
\newcommand{\Ht}{\mathcal H}
\newcommand{\ind}{\mathbf 1}
\newcommand{\dcom}{d_{\mathrm{com}}}
\newcommand{\norm}[1]{\lVert #1\rVert_\infty}
\newcommand{\head}[1]{\par\medskip\noindent\textbf{#1.}\enspace}

\begin{document}
\embedfile[stringmethod=escape,mimetype=text/x-python,
 desc={Python enumeration for the five exceptional types}]{exceptional_entry_boxes.py}
\thispagestyle{plain}
\begin{center}
{\LARGE\bfseries Commuting pairs in lattices\\of reductive Lie algebras\par}
\renewcommand{\author}[1]{\medskip{\large #1\par}}
\newcommand{\address}[1]{{\footnotesize #1\par}}
\newcommand{\email}[1]{{\footnotesize\href{mailto:#1}{\nolinkurl{#1}}\par}}

\author{Caihua Luo}
\address{SCHOOL OF SCIENCE AND ENGINEERING,\\
  THE CHINESE UNIVERSITY OF HONG KONG, SHENZHEN,\\
  GUANGDONG 518172, CHINA.}
\email{luocaihua@cuhk.edu.cn}

\author{Xu Song}
\address{SCHOOL OF MATHEMATICAL SCIENCES, HEBEI NORMAL UNIVERSITY,\\
  SHIJIAZHUANG, HEBEI 050024, CHINA.}
\email{songx@hebtu.edu.cn}

\author{Yilin Zeng}
\address{SCHOOL OF SCIENCE AND ENGINEERING,\\
  THE CHINESE UNIVERSITY OF HONG KONG, SHENZHEN,\\
  GUANGDONG 518172, CHINA.}
\email{223015132@link.cuhk.edu.cn}
\end{center}
\medskip\noindent\textbf{\abstractname.}\enspace
We count ordered commuting pairs in a full lattice \(\Lambda\) in a reductive Lie algebra
\(\g\) over \(\Q\). For a fixed convex body \(B\) containing the origin in its interior,
the count in \(HB\) is \(C_{\Lambda,B}H^{\dim\g+\dim Z(\g)}\) with a power-saving error.
The leading term comes from pairs whose projections to each simple ideal are linearly dependent over its
centroid. Their slopes index an absolutely convergent sum of section volumes for \(C_{\Lambda,B}\).
Subtracting the contribution of dependent pairs from the finite-field Fourier sums gives
cancellation at generic dependent frequencies. For independent frequencies, we bound the
dimension of the locus where the corresponding hyperplane contains an irreducible
component of a nonzero bracket fibre.
For \(n\ge3\), the number of ordered commuting pairs in \(M_n(\Z)\) with entries of
absolute value at most \(H\) is
\[
 \frac{2(n+3)\zeta(n^2-2)}{(n+1)\zeta(n^2-1)}
 (2\lfloor H\rfloor+1)^{n^2+1}
 +O_n\!\left(H^{n^2+1-1/(n+1)}\right),
\]
confirming the exponent \(n^2+1\) predicted by Browning, Sawin and Wang~\cite{BSW}.
\par\medskip
\noindent\textbf{Keywords}\quad Commuting pairs \(\cdot\) Lie lattices \(\cdot\)
Exponential sums \(\cdot\) Integral points

\smallskip\noindent\textbf{Mathematics Subject Classification (2020)}\quad
11D45 \(\cdot\) 11H06 \(\cdot\) 17B20 \(\cdot\) 15B36

\section{Introduction}\label{sec:intro}
Let \(\g\) be a reductive Lie algebra over \(\Q\) and let \(\Lambda\) be a full lattice in \(\g\).
Given a compact convex body \(B\subset\g_{\R}\) with the origin in its interior, we wish to estimate
\[
 N_{\Lambda,B}(H)=\#\{(x,y)\in(\Lambda\cap HB)^2:[x,y]=0\}.
\]
The centre \(\z=Z(\g)\) gives a lower bound at once. If \(d=\dim_{\Q}\g\) and
\(z=\dim_{\Q}\z\), then the pairs with \(y\in\z\) alone give
\[
 N_{\Lambda,B}(H)\ge\#(\Lambda\cap HB)\,\#(\Lambda\cap\z\cap HB)\gg_{\Lambda,B} H^{d+z}.
\]
There are other families of the same size. In \(M_m(\Q)\) one may take \(y=tx+cI_m\) for
any fixed rational \(t\). For a general reductive algebra the slope may vary from one
simple ideal to the next.

To describe these slopes, write
\begin{equation}\label{eq:decomp}
 \g=\z\oplus\bigoplus_{i=1}^s\Res_{\K_i/\Q}\s_i.
\end{equation}
Here \(\K_i\) is the centroid of the \(i\)th rational simple ideal, namely the field of
endomorphisms commuting with its adjoint action, and \(\s_i\) is that ideal viewed as an
absolutely simple \(\K_i\)-algebra. For \(\ell_i=[a_i:b_i]\in\PP^1(\K_i)\) the pairs of slope
\(\ell_i\) have the form \((a_iw_i,b_iw_i)\). Thus a tuple \(\ell=(\ell_i)_i\) gives the
rational subspace
\begin{equation}\label{eq:sections}
 V_\ell=\z^2\oplus\bigoplus_i\{(a_iw_i,b_iw_i):w_i\in\s_i\}\subset\g^2 .
\end{equation}
Every pair in \(V_\ell\) commutes. Its two central entries are free, so \(\dim_{\Q}V_\ell=E\),
where \(E=d+z\). Write \(\mu_\ell\) for the measure on \(V_{\ell,\R}\) in which
\(\Lambda^2\cap V_\ell\) has covolume one.
\begin{theorem}\label{th:main}
There is an \(\eta>0\) such that, for \(H\ge2\),
\begin{equation}\label{eq:main}
 N_{\Lambda,B}(H)=C_{\Lambda,B}H^{d+z}+O_{\Lambda,B}(H^{d+z-\eta}).
\end{equation}
The leading coefficient is positive and is given by the absolutely convergent sum
\begin{equation}\label{eq:coefficient}
 C_{\Lambda,B}=\sum_{\ell\in\prod_i\PP^1(\K_i)}\mu_\ell(B^2\cap V_{\ell,\R}).
\end{equation}
The remainder satisfies the explicit bound \eqref{eq:error} below.
\end{theorem}
\begin{remark}
Let \(\h\) be reductive over a number field \(\F\) of degree \(q\), and let
\(\mathcal M\) be a full \(\OO_\F\)-lattice in \(\h\).
Embed \(\mathcal M\) diagonally in \(\prod_{v\mid\infty}(\h\otimes_\F \F_v)\),
and take \(HB\) in this real vector space, where \(B\) is a fixed compact convex body
containing \(0\) in its interior. Applying Theorem~\ref{th:main} to \(\Res_{\F/\Q}\h\) gives
\[
\#\{(x,y)\in(\mathcal M\cap HB)^2:[x,y]=0\}
=C_{\mathcal M,B}H^{q(d_\F+z_\F)}+O_{\mathcal M,B}(H^{q(d_\F+z_\F)-\eta}),
\]
where \(d_\F=\dim_\F\h\), \(z_\F=\dim_\F Z(\h)\) and \(\eta>0\). Indeed restriction of
scalars multiplies both dimensions by \(q\) and makes \(\mathcal M\) a full \(\Z\)-lattice.
For a simple \(\F\)-ideal of centroid \(\K_i\), the slope still lies in \(\PP^1(\K_i)\),
and the degree in \eqref{eq:error} is \([\K_i:\Q]\). For instance, \((x,ix)\) for
\(0\ne x\in\slie_2(\Q(i))\) has slope \([1:i]\) although its two entries are
independent over \(\Q\).
\end{remark}

When \(\g\) is abelian the sum has one term, \(\vol_\Lambda(B)^2\), where
\(\vol_\Lambda\) assigns covolume one to \(\Lambda\). This includes the zero-dimensional
case, in which the count is one. At the other extreme, for an absolutely simple algebra
and a symmetric body all the sections are dilates of \(B\) itself. We state this
consequence in the normalization used for boxes. Put
\[
 L=2\lfloor H\rfloor+1,\qquad K_{\Lambda,B}=2^{-E}C_{\Lambda,B}.
\]
Thus \(K_{\Lambda,B}\) is the coefficient of \(L^E\). For \(E\ge1\) the change from \(H\)
to \(L\) costs \(O(H^{E-1})\). For absolute rank \(r\ge2\) and dimension \(n\) the saving is
\begin{equation}\label{eq:delta}
 \delta(n,r)=\min\left\{\frac1{r+1},\,1-\frac{r(r+1)}n\right\}.
\end{equation}
It equals \(1/(r+2)\) in type \(A_r\) and \(1/(r+1)\) in the other types (with \(D_3=A_3\)).
\begin{corollary}\label{cor:simple}
Let \(\g\) be absolutely simple over \(\Q\), of dimension \(n\) and rank \(r\).
For a centrally symmetric body \(B=-B\) put \(V_{\Lambda,B}=2^{-n}\vol_\Lambda(B)\),
with \(\vol_\Lambda\) normalized so that \(\Lambda\) has covolume one. Then
\begin{equation}\label{eq:simple}
 N_{\Lambda,B}(H)=4V_{\Lambda,B}\frac{\zeta(n-1)}{\zeta(n)}L^n+
 \begin{cases}O_{\Lambda,B}(H^2\log(2H)),&r=1,\\
 O_{\Lambda,B}(H^{n-\delta(n,r)}),&r\ge2.\end{cases}
\end{equation}
\end{corollary}
Indeed a primitive slope of height \(j\) cuts out a copy of \(j^{-1}B\). Summing its
volume gives the zeta quotient, as we verify in Section~\ref{sec:one}.
In rank one the error \(O_{\Lambda,B}(H^2\log(2H))\) holds for every centrally
symmetric body \(B=-B\). For the unit cube in a lattice basis,
Proposition~\ref{pr:cube-error} improves it to \(O(H^2(\log(2H))^{2/3})\),
as recorded in \eqref{eq:basis}.

For a full matrix algebra the scalar coordinate varies along each section and
changes this coefficient.
\begin{corollary}\label{cor:matrix}
Let \(N_n(H)\) count the ordered commuting pairs in \(M_n(\Z)\) with all entries of
absolute value at most \(H\). For \(n\ge3\),
\[
 N_n(H)=\frac{2(n+3)\zeta(n^2-2)}{(n+1)\zeta(n^2-1)}\,L^{n^2+1}+O_n(H^{n^2+1-1/(n+1)}).
\]
For \(n=2\) the same main term holds with error \(O(H^4\log(2H))\).
\end{corollary}
Browning, Sawin and Wang~\cite[after Theorem~1.1]{BSW} predicted that the lower
bound \(N_n(H)\gg_n H^{n^2+1}\) has the correct exponent. In his 2025 Oberwolfach
report, Browning~\cite[p.~2780]{Browning} records the conjecture that
\(N_n(H)\ll_{n,\varepsilon} H^{n^2+1+\varepsilon}\) for every \(\varepsilon>0\).
Chapman and Mudgal~\cite{CM} obtained an asymptotic formula for \(n=2\).
For \(n=3\), they used restricted divisor correlations to prove \(N_3(H)\asymp H^{10}\).
Corollary~\ref{cor:matrix} proves the conjecture for every \(n\ge4\) and strengthens
the \(n=3\) result to an asymptotic formula with a power-saving error.

Browning, Sawin and Wang~\cite[Theorem~1.1]{BSW} obtained
\(N_n(H)\ll_n H^{n^2+2-2/(n+1)}\).
For \(n=2\) our coefficient \(\tfrac{10}{3}\zeta(2)/\zeta(3)\) and error term agree
with those of Chapman and Mudgal. For \(n=3\) the coefficient is
\(3\zeta(7)/\zeta(8)\), and the error is \(O(H^{39/4})\).
The lower bound of Chapman and Mudgal in every size includes the family \(y=x+cI_n\).
Corollary~\ref{cor:matrix} accounts for every rational slope modulo the scalars.
Mudgal~\cite{Mudgal} also treats commuting pairs in arbitrary finite sets of real
\(2\times2\) matrices.

There is a corresponding comparison for trace-zero matrices. Adding independent
integral scalars to the two matrices embeds \(L^2\) disjoint copies of the trace-zero
count at height \(H\) in the full count at height \(2H\). The traces determine both
shifts. Hence the bound of~\cite{BSW} implies \(O_n(H^{n^2-2/(n+1)})\) in the
trace-zero case. Corollary~\ref{cor:simple} gives the leading term of order \(H^{n^2-1}\),
with a coefficient computed in Appendix~\ref{sec:boxes}.

\head{Outline of the proof}
Let \(x_i,y_i\) be the projections of \(x,y\) to \(\s_i\). Then
\[
 (x,y)\in\bigcup_\ell V_\ell
 \quad\Longleftrightarrow\quad
 x_i\wedge y_i=0\ \text{in }\bigwedge\nolimits^2_{\K_i}\s_i
 \quad(1\le i\le s).
\]
A pair in this union has a unique slope tuple if every \((x_i,y_i)\ne(0,0)\).
Summing the uniform lattice estimate over the slopes and counting zero projections
separately gives \eqref{eq:coefficient}. Outside the union it
suffices to save a power of \(H\) in one factor with \(x_i\wedge y_i\ne0\).

We complete the count modulo a good split prime \(p\) as in~\cite[\S4]{BSW}. In a factor of
dimension \(n\), the full transform \(S\) and the transform \(S^\circ\) restricted
to independent pairs satisfy, at a generic nonzero dependent frequency,
\[
 S=p^n+O(p^{n-1/2}),\qquad S^\circ=S-(p^n-p)=O(p^{n-1/2}).
\]
At independent frequencies, a hyperplane may contain a component of the bracket fibre.
Proposition~\ref{pr:geometry}(ii) bounds the dimension of this exceptional set, and
averaging over primes bounds the number of short frequencies reducing into it.

Sections~\ref{sec:coeff} and~\ref{sec:slopes} treat the slope sections.
Section~\ref{sec:fourier} estimates the independent pairs using the bracket geometry
of Section~\ref{sec:geom}. Section~\ref{sec:remainder} treats the error term, and
Appendix~\ref{sec:boxes} evaluates the box coefficients.

\section{The leading coefficient}\label{sec:coeff}
We evaluate \eqref{eq:coefficient} for symmetric bodies. We write \(C_{\Lambda,B}\)
for the coefficient of \(H^E\) and \(K_{\Lambda,B}=2^{-E}C_{\Lambda,B}\) for the
coefficient of \(L^E\). Most of the examples reduce to dilating a simple factor or
taking a slice along the centre.

\subsection{One simple factor and products}\label{sec:one}
\begin{proof}[Proof of Corollary~\ref{cor:simple}]
Represent a rational slope by coprime integers \(a,b\), up to a common sign, and put
\(j=\max(|a|,|b|)\). If \(aw,bw\in\Lambda\), Bézout's identity gives \(w\in\Lambda\).
The two conditions \(aw,bw\in B\) reduce to \(w\in j^{-1}B\) by symmetry and convexity.
Thus a slope of height \(j\) contributes \(j^{-n}\vol_\Lambda(B)\) to
\eqref{eq:coefficient}. There are \(4\varphi(j)\) primitive slopes of height \(j\),
including the four slopes of height one. Consequently
\[
 C_{\Lambda,B}=4\vol_\Lambda(B)\sum_{j\ge1}\frac{\varphi(j)}{j^n}
 =4\vol_\Lambda(B)\frac{\zeta(n-1)}{\zeta(n)}.
\]
Changing the normalization and applying \eqref{eq:error} gives \eqref{eq:simple}.
The argument applies to nonsplit absolutely simple forms over \(\Q\) as well.
\end{proof}
Suppose next that both the lattice and the box are products in \eqref{eq:decomp},
with the cube of a lattice basis on each summand. The count then factors exactly.
The central coordinates contribute \(L^{2z}\), and each simple factor contributes
its own commuting-pair count. When all the centroid fields are \(\Q\),
Corollary~\ref{cor:simple} therefore gives
\begin{equation}\label{eq:product}
 K_{\Lambda,B}=\prod_{i=1}^s\frac{4\zeta(n_i-1)}{\zeta(n_i)},\qquad n_i=\dim_{\K_i}\s_i.
\end{equation}
For an abelian algebra the empty product is one.

\subsection{A simple quotient and its central slices}
Now suppose that \(\g/\z\) is absolutely simple over \(\Q\), of dimension \(m\).
Represent a slope \([a:b]\) by coprime integers. Interchange the entries if necessary
so that \(j=|a|=\max(|a|,|b|)\), and put \(t=|b|/j\). For fixed \(x\in B\), the possible
second entries on this section are translates of \((b/a)x\) by elements of \(\z_\R\).
By symmetry their volume is \(\int_{\z_\R}\ind_B(tx+c)\,dc\). This leads to the slice integral
\begin{equation}\label{eq:slice}
 f_B(t)=2^{-E}\int_B\int_{\z_\R}\ind_B(tx+c)\,dc\,dx,\qquad 0\le t\le1,
\end{equation}
where \(dx\) and \(dc\) give covolume one to \(\Lambda\) and \(\Lambda\cap\z\), respectively.

The contribution of this slope to \(K_{\Lambda,B}\) is \(j^{-m}f_B(t)\). To see the
factor \(j^{-m}\), project the section lattice onto its first entry. In the free quotient
\(\bar\Lambda=\Lambda/(\Lambda\cap\z)\) the equation is \(a\bar y=b\bar x\).
It has a solution precisely when \(\bar x\in a\bar\Lambda\), since \((a,b)=1\).
The image consequently has index \(j^m\) in \(\Lambda\). The kernel is the central
lattice in the second coordinate. Integrating first over this kernel gives the assertion.

For \(j>1\), pairing a reduced residue \(u\) with \(j-u\) shows that the mean of \(t\)
at fixed height is \(1/2\). At \(j=1\) the four ratios are \(0,0,1,1\). Hence an affine
slice function gives the closed formula
\begin{equation}\label{eq:affine}
 K_{\Lambda,B}=4f_B(1/2)\frac{\zeta(m-1)}{\zeta(m)}\qquad\text{if \(f_B\) is affine.}
\end{equation}
For general \(f_B\) one can sum over all nonzero pairs modulo sign and divide by
\(\sum_{j\ge1}j^{-m}\) to remove common multiples. The two endpoints each occur twice
on the boundary of the square of radius \(j\), and each interior value \(u/j\) occurs
four times. This gives
\begin{equation}\label{eq:allslices}
 K_{\Lambda,B}=\frac{4}{\zeta(m)}\sum_{j\ge1}j^{-m}
 \left(\frac{f_B(0)+f_B(1)}2+\sum_{u=1}^{j-1}f_B(u/j)\right).
\end{equation}
\begin{proof}[Proof of Corollary~\ref{cor:matrix}]
For \(x\) uniform in the entry box \(B\) of \(M_n(\R)\) and \(0\le t\le1\),
\[
 tx+cI_n\in B
 \quad\Longleftrightarrow\quad
 -1-t\min_i x_{ii}\le c\le1-t\max_i x_{ii}.
\]
Set \(R(x)=\max_i x_{ii}-\min_i x_{ii}\) and \(U_i=(1+x_{ii})/2\).
The \(U_i\) are independent and uniform on \([0,1]\), so
\[
 \mathbb E R(x)=2\int_0^1\bigl(1-u^n-(1-u)^n\bigr)\,du
 =\frac{2(n-1)}{n+1}.
\]
The interval for \(c\) has length \(2-tR(x)\). Thus \eqref{eq:slice} gives
\(f_B(t)=1-\tfrac12t\,\mathbb E R(x)\), and \eqref{eq:affine} yields
\[
 K_{\Lambda,B}=4f_B(1/2)\frac{\zeta(n^2-2)}{\zeta(n^2-1)}=\frac{2(n+3)\zeta(n^2-2)}{(n+1)\zeta(n^2-1)}.
\]
The error follows from \eqref{eq:error} for the quotient of type \(A_{n-1}\).
\end{proof}
The first shell already explains the lower bound of Chapman and Mudgal
\cite[Lemma~1.4]{CM}. Its slopes \(0\) and \(\infty\) contribute \(1\) each, and \(1\)
contributes \(f_B(1)=2/(n+1)\). The slope \(-1\) contributes a further \(2/(n+1)\).
The shell of height \(j\) contributes \(\varphi(j)j^{-(n^2-1)}\) times the first shell,
accounting for \(\zeta(n^2-2)/\zeta(n^2-1)\) in Corollary~\ref{cor:matrix}.

\subsection{A formula in arbitrary matrix coordinates}
Even when a lattice couples different ideals, its section volumes can be computed in
the original matrix entries. Let \(A\) have as its columns the entries of the basis
matrices of \(\Lambda\), in a fixed order. In these coordinates the entry box is
\(\{u:\norm{Au}\le1\}\). Write the slope equations \(b_ix_i-a_iy_i=0\) in rational
coordinates and clear denominators. If the columns of \(U_\ell\) form a \(\Z\)-basis
of the full integer kernel of these equations in the coordinates of \(\Lambda^2\), then
\begin{equation}\label{eq:kernel}
 K_{\Lambda,B}=2^{-E}\sum_\ell\vol_E\{v\in\R^E:\norm{\operatorname{diag}(A,A)U_\ell v}\le1\}.
\end{equation}
Indeed \(v\mapsto U_\ell v\) carries \(\Z^E\) onto the section lattice, so ordinary
volume in \(v\) is precisely \(\mu_\ell\). Smith normal form computes this kernel,
including any congruences between ideals. Each summand is the volume of a rational
polytope. Truncation at \(\max_i\Ht(\ell_i)\le V\) costs
\(O_{\Lambda,B}(\sum_i V^{-e_i(n_i-2)})\) by \eqref{eq:tail}, where \(e_i=[\K_i:\Q]\)
and \(\Ht\) is the projective height defined in Section~\ref{sec:slopes}.

\section{Counting pairs of fixed slope}\label{sec:slopes}
Let \(e_i=[\K_i:\Q]\), \(n_i=\dim_{\K_i}\s_i\), and let \(r_i\) be the absolute rank of
\(\s_i\). Thus \(E=2z+\sum_i e_i n_i\). In this notation the quantitative assertion of
Theorem~\ref{th:main} is
\begin{equation}\label{eq:error}
 \big|N_{\Lambda,B}(H)-C_{\Lambda,B}H^E\big|\ll_{\Lambda,B}
 H^{E-1}\mathcal L(H)+\sum_{i:r_i\ge2}H^{E-e_i\delta(n_i,r_i)},
\end{equation}
where \(\mathcal L(H)=\log(2H)\) if some \((e_i,n_i)=(1,3)\) and is one otherwise.
The first term comes from counting lattice points on the slope sections. Each term in
the sum corresponds to pairs that are independent in the \(i\)th simple factor.

We recall why the centroid gives the slope field in \eqref{eq:decomp}. For a rational
simple ideal its centroid is a division algebra by simplicity and is commutative by
perfectness. The centroid is defined by linear equations, so its formation commutes
with extending scalars. Over \(\overline{\Q}\) it is the product of one scalar field
for each geometric simple ideal. Its primitive idempotents therefore identify those
ideals with the embeddings of the centroid, giving absolute simplicity over that field.

It suffices to estimate independent pairs in a single simple factor. For a number field
\(\K\) of degree \(e\), use the usual absolute value at every embedding, counting conjugate
complex embeddings twice in products. In fixed \(\K\)-coordinates put
\(\norm{x}=\max_{\sigma,j}|\sigma x_j|\).

\begin{proposition}\label{pr:independent}
Let \(\K\) have degree \(e\) and let \(\s/\K\) be absolutely simple, of dimension \(n\)
and absolute rank \(r\). In fixed \(\K\)-coordinates the number of commuting,
\(\K\)-linearly independent pairs in \((\OO_\K^n)^2\) with both norms at most \(H\) is
\[
 \ll_{\s,\K} H^{en-e\delta(n,r)}\qquad(r\ge2).
\]
For \(r=1\) this number is zero.
\end{proposition}
In rank one, \(\bigwedge_\K^2\s\to\s\) is surjective by perfectness and has source
and target of dimension three. The vanishing of the bracket is thus equivalent to
dependence. We prove the higher-rank assertion in Section~\ref{sec:fourier}.

\subsection{Heights and integral representatives}\label{sec:heights}
For \(\ell=[a:b]\in\PP^1(\K)\) use the absolute multiplicative height
\[
 \Ht(\ell)^e=\frac{\prod_{\sigma:\K\hookrightarrow\C}\max(|\sigma a|,|\sigma b|)}
 {\mathrm N(a\OO_\K+b\OO_\K)} .
\]
The denominator is the absolute ideal norm. Choose one integral ideal in each ideal
class and scale \((a,b)\) until \(\id=a\OO_\K+b\OO_\K\) is one of these representatives.
We still have the freedom to multiply both entries by a unit. In logarithmic coordinates
the unit lattice has bounded fundamental domain in the hyperplane of weighted mean
zero. We can therefore keep each local maximum within a bounded factor of its
geometric mean. Since the possible norms of \(\id\) are fixed, this gives
\begin{equation}\label{eq:balanced}
 \max(|\sigma a|,|\sigma b|)\asymp_\K\Ht(\ell)\qquad\text{at every embedding \(\sigma\).}
\end{equation}
Counting representatives in two copies of the Minkowski lattice gives
\begin{equation}\label{eq:heightcount}
 \#\{\ell:\Ht(\ell)\le V\}\ll_\K V^{2e}\qquad(V\ge1).
\end{equation}
The same balancing argument applies to a projective vector with any fixed number of coordinates.

For a balanced representative, a nonzero integral pair of slope \(\ell\) is uniquely
\((aw,bw)\) with \(0\ne w\in(\id^{-1})^n\). Above the box of radius \(H\) we have
\(\norm{w}\ll H/\Ht(\ell)\). The finite list of parameter lattices has a positive
minimum length. Hence only slopes with \(\Ht(\ell)\ll H\) contribute, and their
individual contributions are \(O((H/\Ht(\ell))^{en})\). For \(V\ge1\), divide the
remaining heights into shells \(R<\Ht(\ell)\le2R\), with \(R=2^jV\). By
\eqref{eq:heightcount} their combined contribution is at most
\[
 \sum_{\substack{R=2^jV\ll H\\j\ge0}}R^{2e}(H/R)^{en}\ll H^{en}V^{-e(n-2)}.
\]
Thus all dependent pairs together cost \(O(H^{en})\), including the zero pair and
the bounded-height slopes. Fixed fractional lattices satisfy the same estimates after
clearing denominators.

\subsection{The lattice discrepancy}
Fix \(\ell\), write \(R_i=\Ht(\ell_i)\), \(d_i=e_i n_i\) and put
\(\lambda_\ell=\mu_\ell(B^2\cap V_{\ell,\R})\). In the coordinates
\((u,v,(w_i))\) of \eqref{eq:sections} the natural comparison lattice
\(\Gamma=\Z^z\oplus\bigoplus_i\OO_{\K_i}^{n_i}\) has section lattice
\[
 \Gamma_\ell=\Z^{2z}\oplus\bigoplus_i(\id_i^{-1})^{n_i}.
\]
Choose \(q\) so that \(q\Gamma\subset\Lambda\subset q^{-1}\Gamma\). The pullback of
\(\Lambda^2\cap V_\ell\) lies between \(q\Gamma_\ell\) and \(q^{-1}\Gamma_\ell\).
The ideals \(\id_i\) lie on fixed lists, and each quotient
\(q^{-1}\Gamma_\ell/q\Gamma_\ell\) is finite. Only finitely many section lattices occur in
these coordinates. We may therefore choose fundamental cells of uniformly bounded
diameter and covolume bounded away from zero.

In real Minkowski coordinates let \(P_\ell(H)\) be the pullback of \((HB)^2\). It lies
between two fixed positive multiples of the product body with central radii \(H\) and
respective factor radii \(H/R_i\). Indeed \(B\) lies between fixed multiples of a
product of coordinate balls. For a complex coordinate the maximum of \(|a_iw|,|b_iw|\)
is \(\max(|a_i|,|b_i|)|w|\), and \eqref{eq:balanced} applies at each place. In particular
\begin{equation}\label{eq:tail}
 0<\lambda_\ell\ll\prod_i R_i^{-d_i}.
\end{equation}
Equations \eqref{eq:heightcount} and \eqref{eq:tail} prove convergence and positivity
in \eqref{eq:coefficient}, and bound its tail over \(\max_i R_i>V\) by
\(O(\sum_i V^{-e_i(n_i-2)})\) when \(s>0\).

Choose \(c_0\) sufficiently large that any pair with every \(w_i\ne0\) has \(R_i\le c_0H\)
for all \(i\). In this range put \(R_*=\max(1,R_1,\ldots,R_s)\). The number of points
on this slope with every \(w_i\ne0\) equals
\begin{equation}\label{eq:discrepancy}
 H^E\lambda_\ell+O\left(H^{E-1}R_*\prod_iR_i^{-d_i}\right).
\end{equation}
Put \(P=P_\ell(H)\) and let \(\mathcal U\) be the union of fundamental cells centred
at lattice points of \(P\). Each cell lies within distance \(\rho=O(1)\) of its
centre, and \(P\) contains a ball of radius \(s_0\gg H/R_*\) about zero. Convexity gives
\[
 (1-t)P\subseteq\mathcal U\subseteq(1+t)P
 \qquad\bigl(t=\rho/s_0<1\bigr).
\]
The upper inclusion holds for all \(t\ge0\). Since \(t\ll R_*/H\ll1\), the volume
discrepancy is \(O\bigl((R_*/H)H^E\prod_iR_i^{-d_i}\bigr)\).
Removing the points with \(w_i=0\) costs at most \(O(H^{2z}\prod_{j\ne i}(H/R_j)^{d_j})\)
for each \(i\), by the product lattice comparison. As \(R_i/H\le c_0\) and \(d_i\ge1\),
this also fits the error in \eqref{eq:discrepancy}.

\begin{proof}[Proof of Theorem~\ref{th:main}, assuming Proposition~\ref{pr:independent}]
Write \(N_{\rm dep}(H)\) for the pairs dependent in every simple factor. First count
those with every projected pair nonzero, so that each pair has exactly one slope tuple.
To sum the error in \eqref{eq:discrepancy}, use \(R_*\le1+\sum_iR_i\). All factors
except the one carrying the extra \(R_i\) have a convergent height sum. In that factor
\eqref{eq:heightcount} gives
\[
 \sum_{\substack{\ell_i\\R_i\le c_0H}}R_i^{1-d_i}
 \ll\sum_{0\le j\le\log_2(2c_0H)}2^{j(1-e_i(n_i-2))}\ll\mathcal L(H).
\]
There is a logarithm only for \(e_i=1,\ n_i=3\). Extending the main sum to all \(\ell\)
costs \(O(H^{E-1})\) by \eqref{eq:tail} and its tail estimate. If the projection of the pair to
the \(i\)th factor vanishes, the other dependent pairs and the centre give
\(O(H^{E-d_i})\) possibilities. Consequently
\[
 N_{\rm dep}(H)=H^E\sum_\ell\lambda_\ell+O_{\Lambda,B}(H^{E-1}\mathcal L(H)).
\]
Now fix an index \(i\) at which the pair is independent. Enlarge the lattice and the
body to the containing product lattice and box used above. In every other factor
the slope estimate and Proposition~\ref{pr:independent} together give \(O(H^{d_j})\)
choices. In rank one there are no independent commuting pairs. Summing over the
possible indices therefore gives
\[
 0\le N_{\Lambda,B}(H)-N_{\rm dep}(H)\ll H^{2z}\sum_{i:r_i\ge2}
 H^{d_i-e_i\delta(n_i,r_i)}\prod_{j\ne i}H^{d_j}
 =\sum_{i:r_i\ge2}H^{E-e_i\delta(n_i,r_i)}.
\]
This is \eqref{eq:error}. When there are no simple factors the ordinary lattice count
in \((HB)^2\) gives the assertion directly.
\end{proof}

\section{Fourier estimates for independent pairs}\label{sec:fourier}
We turn to Proposition~\ref{pr:independent}. Fix an absolutely simple \(\K\)-algebra
\(\s\) of dimension \(n\) and rank \(r\ge2\), with \([\K:\Q]=e\). Put
\[
 M=2n,\qquad \dcom=n+r,\qquad k=n-r .
\]
Here \(M\) is the dimension of the space of pairs and \(\dcom\) the dimension of the
commuting variety. We discard the finitely many primes dividing denominators of the
bracket or at which the Killing form \(\kappa\) degenerates. We first compute the
complete Fourier sum and then estimate the frequencies that occur after Poisson summation.

\subsection{Finite-field transforms}
For \(a,b\) in the Lie algebra over the field under consideration, put
\[
 F_a=\{(y,u):[y,u]=a\},\qquad F_{a,b}=\{(y,u)\in F_a:\kappa(b,y)=0\}.
\]
Work first with a split model \(\s_p\) over a good prime field \(\F_p\). Fix a
nontrivial additive character \(\psi\) and put
\[
 S(a,b)=\sum_{\substack{x,y\in\s_p\\{[x,y]=0}}}\psi\big(\kappa(a,x)+\kappa(b,y)\big).
\]
Write \(S^\circ\) for the same sum restricted to linearly independent pairs. Insert
the additive-character indicator for \([x,y]=0\) and sum first over \(x\).
Invariance and nondegeneracy of the Killing form give
\begin{equation}\label{eq:transform}
 \begin{aligned}
 S(a,b)&=p^{-n}\sum_{y,u}\psi(\kappa(b,y))\sum_x\psi\big(\kappa(x,a+[y,u])\big)\\
 &=\sum_{(y,u)\in F_{-a}(\F_p)}\psi(\kappa(b,y)).
 \end{aligned}
\end{equation}
For \(a\ne0\), average \((y,u)\mapsto(ty,t^{-1}u)\) over \(t\in\F_p^\times\).
The identity \(\sum_{t\in\F_p^\times}\psi(tc)=p\ind_{\{c=0\}}-1\), applied with
\(c=\kappa(b,y)\), gives
\begin{equation}\label{eq:average}
 S(a,b)=\frac{p\,\#F_{-a,b}(\F_p)-\#F_{-a}(\F_p)}{p-1}.
\end{equation}
We next compute the term that has to be subtracted. The dependent cone \(\mathcal D\)
is the union of its \(p+1\) slope subspaces. Subtract \(p\) copies of the origin to
correct their common intersection. Summing first along each subspace gives the exact identity
\[
 \widehat{\ind_{\mathcal D}}(a,b)=p^n\#\{[u:v]\in\PP^1(\F_p):ua+vb=0\}-p .
\]
For \((a,b)\ne(0,0)\) the projective count is one or zero according as \(a,b\) are
dependent or independent. In the dependent case a change of the two variables by
\(\GL_2(\F_p)\) takes the frequency to \((a_0,0)\) on the same line. Thus
\eqref{eq:transform} gives, for \((a,b)\ne(0,0)\),
\[
 S^\circ(a,b)=\begin{cases}
 \#F_{-a_0}(\F_p)-p^n+p,&a,b\text{ dependent},\\
 S(a,b)+p,&a,b\text{ independent}.
 \end{cases}
\]
For a dependent frequency the leading \(p^n\) cancels as soon as \(F_{-a_0}\) is
geometrically irreducible of dimension \(n\). For an independent frequency,
\eqref{eq:average} requires an upper bound for the hyperplane section. The following
proposition supplies both geometric results. We prove it in Section~\ref{sec:geom}.

\begin{proposition}\label{pr:geometry}
Let \(\s\) be simple of dimension \(n\) and rank \(r\) over an algebraically closed
field of characteristic zero. For \(a\ne0\) the variety \(F_a\) is nonempty and pure
of dimension \(n\). The first coordinate is regular semisimple on a dense open subset
of each component. The following bounds hold.
\begin{enumerate}[label=\textup{(\roman*)}]
\item Outside a proper closed cone \(J\subset\s\), the fibre \(F_a\) is geometrically irreducible.
\item The closure of the independent pairs \((a,b)\) with \(\dim F_{a,b}=n\) has dimension at most \(n\).
\end{enumerate}
At the origin, \(\dim F_0=n+r\).
\end{proposition}
Browning, Sawin and Wang proved that the complex matrix commutator map
\(M_m(\C)^2\to\slie_m(\C)\) is flat away from zero \cite[Theorem~1.2]{BSW}.
Writing \(M_m(\C)=\slie_m(\C)\oplus\C I_m\) identifies its fibre over \(a\)
with \(F_a\times\A^2\), where \(F_a\) is the bracket fibre in \(\slie_m(\C)^2\).
For \(a\ne0\), the dimension assertion therefore gives
\((m^2-1)+2=m^2+1\) for the matrix fibre, as in their theorem.

Let \(Y\) be the closed exceptional set in part (ii), with \(a\) replaced by \(-a\) as
in \eqref{eq:transform}. Choose models of \(J\) and \(Y\) over a localization of the
ring of integers of a field of definition.

\begin{lemma}\label{le:finite}
Let \(\p\) be a prime at which these models have good reduction, with residue field
\(\F_p\), and suppose that the reduction of \(\s\) at \(\p\) is split. For \(p\)
sufficiently large, the following bounds hold uniformly in the nontrivial character
\(\psi\).
\begin{equation}\label{eq:crude}
 |S^\circ(0,0)|\ll p^{\dcom},\qquad |S^\circ(a,b)|\ll p^n\quad((a,b)\ne(0,0)).
\end{equation}
For a nonzero dependent frequency with common line outside \(J\),
\begin{equation}\label{eq:dep}
 |S^\circ(a,b)|\ll p^{n-1/2}.
\end{equation}
For independent \(a,b\) with \((a,b)\notin Y\),
\begin{equation}\label{eq:ind}
 |S(a,b)|\ll p^{n-1}.
\end{equation}
The last bound holds also for \(S^\circ\).
\end{lemma}
The primes chosen using Lemma~\ref{le:cones} have this residue field and split reduction.
\begin{proof}
We first pass from characteristic zero to good reduction. For the bracket and its
hyperplane incidence, the loci where a fibre exceeds the dimensions in
Proposition~\ref{pr:geometry} are constructible. Outside \(J\) the same holds for
failure of geometric irreducibility or of dimension \(n\)~\cite{Stacks}. In each
case impose the stated conditions on \(a,b\), treating \(a=0\) separately.
The resulting bad locus has empty generic fibre. By Chevalley's theorem its image
in the one-dimensional base is constructible and misses the generic point, hence
consists of finitely many primes. Away from these primes all the required geometric
assertions hold simultaneously. The degrees are bounded by the fixed equations.

At zero we count at most \(\#F_0(\F_p)\ll p^{n+r}\). For \(a\ne0\), \eqref{eq:transform}
gives \(|S(a,b)|\le\#F_{-a}(\F_p)\ll p^n\), and symmetry permits us to exchange
\(a,b\). Removing the dependent cone therefore gives \eqref{eq:crude}. Outside \(J\),
uniform Lang--Weil \cite[Theorem~7.7.1]{Poonen} gives
\(\#F_{-a_0}(\F_p)=p^n+O(p^{n-1/2})\). Subtracting \(p^n-p\) proves \eqref{eq:dep}.
For linearly independent \(a,b\) with \((a,b)\notin Y\), the dimension of
\(F_{-a,b}\) is at most \(n-1\). Substitution in \eqref{eq:average} gives
\[
 |S(a,b)|\ll(p\cdot p^{n-1}+p^n)/(p-1)\ll p^{n-1}.
\]
Here \(S^\circ=S+p\), so the same estimate holds for the restricted sum.
\end{proof}

\subsection{Stratification and the choice of primes}\label{sec:strat}
We use the following specialization of Fouvry--Katz \cite[Theorem~1.1]{FK}.
Let \(m,d_V\ge1\), let \(V\subset\A^m_{\Z}\) be a locally closed subscheme with
\(\dim V_{\C}\le d_V\), and fix \(f\in\Z[t_1,\ldots,t_m]\). With multiplicative
weight one (take the unit function \(g=1\) and the trivial character), there are
nested closed subschemes \(X_1\supseteq\cdots\supseteq X_m\) of the dual affine
space, of relative dimension at most \(m-i\) for \(X_i\), and a constant \(A\)
depending only on \(m,d_V,V,f\). For \(1\le i\le m\), a nontrivial additive character
\(\psi\) of \(\F_p\), and \(\xi\notin X_i(\F_p)\), their estimate is
\[
 \left|\sum_{v\in V(\F_p)}\psi\big(f(v)+\xi\cdot v\big)\right|
 \le A p^{(d_V+i-1)/2}.
\]

Take a split Chevalley model and the closed subscheme \(\mathcal C\subset\A^M_\Z\)
defined by the integral bracket equations \([x,y]=0\). Its complex fibre has
underlying variety of dimension \(\dcom=n+r\) by Richardson~\cite{Richardson}.
Thus \(m=M\), \(d_V=\dcom\), and the fixed polynomial is \(f=0\). Every map
from \(\operatorname{Spec}\F_p\) to \(\mathcal C_{\F_p}\) factors uniquely through
its reduction, so nilpotents affect neither the point count nor the sum. For all
\(p\) outside a fixed finite set, the Gram matrix of \(\kappa\) carries the
linear phase to the coordinate pairing by an invertible change of frequency.
Pull back the strata by this change and put \(X_0=\A^M\). For \(0\le j<M\),
the estimate with \(i=j+1\) gives
\begin{equation}\label{eq:FK}
 \dim\big((X_j)_{\overline{\F}_p}\big)\le M-j,\qquad
 |S(\xi)|\le A p^{(\dcom+j)/2}
 \quad\big(\xi\notin X_{j+1}(\F_p),\ p\nmid N\big)
\end{equation}
for one fixed integer \(N\ge1\) and every nontrivial \(\psi\). This is also the
indexing used in \cite[\S4, following (4.2)]{BSW}.

For averaging, a frequency will be a vector over \(\K\) whose reduction we test at
many primes. We need one family of cones defined over \(\K\) to make all these tests.
We also need its components to span the ambient space in order to apply the
short-vector estimate below.

\begin{lemma}\label{le:cones}
For a \(\K\)-form of \(\s\), the proper direction cone \(J\) may be chosen over \(\K\).
The exceptional sets \(Y\) and \(X_j\) may be replaced by closed cones \(Y'\) and
\(X'_j\) defined over \(\K\), with \(\dim Y'\le n\), \(\dim X'_j\le M-j\), and
\(X'_{j+1}\subseteq X'_j\) for \(0\le j<M\). Every positive-dimensional geometric
component of a replacement cone spans the full frequency space.

There is a finite extension \(\K'/\K\) such that, at primes above all sufficiently
large rational primes \(p\) splitting completely in the Galois closure of \(\K'/\Q\),
the residue field is \(\F_p\) and \(\s\) has split reduction. At these primes
\eqref{eq:ind} holds for independent pairs outside the reduction of \(Y'\), and
\eqref{eq:FK} holds with the \(X'_j\) in place of the \(X_j\).
\end{lemma}
\begin{proof}
Over a splitting field the connected group \(G\times\GL_2\) acts on the frequency
space by the dual action. The transform \(S\) is invariant, since the primal action
permutes commuting pairs. For a stratum \(X\) intersect all its translates over
the algebraic closure. By Noetherianity the intersection equals an intersection
of finitely many translates, say \(\bigcap_\nu g_\nu X\). Make these choices for
all \(X_j\) and for \(Y\) at once. For primes splitting in a field of definition
of all \(g_\nu\), a point outside the reduced intersection has a translate outside
the original \(X\). The estimate for the latter point applies to the former by
invariance. For \(Y\) this reasoning is restricted to independent pairs, an invariant
condition. Intersecting all translates preserves the nesting of the \(X_j\).

These intersections are invariant in characteristic zero. Scalar matrices in
\(\GL_2\) make them cones. The connected group fixes each geometric component, whose
linear span is consequently an invariant subspace. The ambient representation is
absolutely irreducible, being the external tensor product of the adjoint representation
and a two-dimensional standard or dual representation. Hence any component containing
a nonzero point spans the ambient space.

The cone \(J\) may be defined as the closure of the locus of nonzero targets where
geometric irreducibility fails, with the origin adjoined. This definition is invariant
under Lie algebra automorphisms and scalars. It gives the same proper direction cone
after transport to the \(\K\)-form, and \eqref{eq:dep} holds off its reduction after
the finite exclusions above. Transport the frequency cones to the chosen \(\K\)-form
over a finite splitting field and take, for each cone, the union of its Galois
conjugates. These unions give the \(\K\)-defined cones \(Y'\) and \(X'_j\).
They preserve the nesting, the dimension bounds and the estimates on the
complements, and each surviving component still spans. All translates, intersections
and conjugates used here are finite in number. Choosing equations over a common
localization and removing finitely many primes preserves the set-theoretic identities
under reduction. Enlarge \(\K'\) to include the splitting field and these fields of
definition. The Gram matrix of \(\kappa\) identifies the frequency coordinates through
the same fixed invertible linear map on each copy of \(\s\) and preserves dependence.
\end{proof}
From now on \(Y\) and \(X_j\) denote the replacement cones of Lemma~\ref{le:cones}.
Choose a finite Galois extension of \(\Q\) containing the fields just used and the
Hilbert class field of \(\K\). Take the rational primes splitting completely in it,
outside the finite exceptional set. Chebotarev's theorem gives \(\gg P/\log P\)
such primes in \([P,2P]\). For each \(p\) choose one prime \(\p\) of \(\K\) above it.
It has residue field \(\F_p\) and is principal by the Hilbert class field criterion.
Write \(\mathscr P(P)\) for the resulting set of prime ideals and \(\Av_\p\) for
the mean over it. The estimates \eqref{eq:crude}--\eqref{eq:FK} hold for every
nontrivial additive character at these primes.

\subsection{Short vectors and reduction modulo primes}\label{sec:short}
We need an averaged bound for \(Y\) and for the two deepest strata in \eqref{eq:FK}
used below. The elementary bound in the first assertion of Lemma~\ref{le:short}
also applies to \(J\subset\A^n\) with \(b=n-1\), and to \(X_j\subset\A^M\) with
\(b=M-j\) for \(0\le j\le k-4\). We state the lemma with an arbitrary ambient
dimension \(m\) to cover both uses. A cone may be reducible. Norms on a fixed
fractional ideal can be used after clearing a fixed denominator.

\begin{lemma}\label{le:short}
Let \(X\subset\A^m_\K\) be closed of dimension at most \(b\). For \(U\ge1\),
put \(W=U^e\). The number of \(\xi\in\OO_\K^m\) with \(\norm{\xi}\le U\) and
\(\xi\bmod\p\in X(\F_p)\) is \(O(W^b)\) if \((2U)^e<p\), at good primes.
Suppose in addition that \(X\) is a cone, \(b\ge3\), every geometric component of
dimension \(b\) spans \(\A^m\), and \(W^2\log P\ll P\). For sufficiently large \(P\),
the mean of this count over the chosen primes is
\[
 \begin{cases}O(W^{b-1}),&m-b\ge4,\\ O(W^{b-1/2}),&m-b=3.\end{cases}
\]
\end{lemma}
\begin{proof}
A finite linear Noether normalization defined over \(\K\) puts the coordinates of all
projected points in a fixed fractional lattice of radius \(O(U)\). There are
\(O(W^b)\) possibilities, with bounded fibres modulo all good primes. Reduction
is injective on the original box, since a nonzero difference in \(\p\) would have
absolute norm at least \(p\) and at most \((2U)^e\). This proves the first assertion
and also handles components of smaller dimension in the average.

For a \(\K\)-irreducible component of dimension \(b\) take its reduced projectivization
and a general \(\K\)-linear projection to \(\PP^b\). The projection has centre disjoint
from the projectivization and is finite and generically one-to-one on its image.
Indeed for a general point on each geometric component the union of secants from
that point to the whole projectivization has closure of dimension at most \(b\).
A general centre of dimension \(m-b-2\) misses this secant locus and the tangent
spaces at the chosen smooth general points. The same incidence argument for
projection from a point is given in \cite[\S43.23, Lemmas 43.23.1--2]{Stacks}.
These are finitely many nonempty open conditions on the Grassmannian, and \(\K\)
is infinite, so the centre can be chosen over \(\K\). Disjointness from the variety
gives finiteness, and the separated general points give generic degree one in
characteristic zero. Taking secants to the whole union ensures that distinct
geometric components have distinct images. Pulling back general hyperplanes shows
that each image component has its original degree, at least \(m-b+1\) by the
degree bound for a nondegenerate projective variety. The induced map of cones
into \(\A^{b+1}\) is finite, with uniformly bounded fibres after good reduction.
Clear a fixed denominator in the projected coordinates, and let \(F\) be an
integral homogeneous equation for the image.

Write \(N_X(U;\p)\) for the count on this component. A projected vector \(w\) either
satisfies \(F(w)=0\) or costs a prime divisor of its nonzero norm
\(|\mathrm N_{\K/\Q}F(w)|\). Bounded fibres and injectivity on the original box give
\[
 \Av_\p N_X(U;\p)\ll
 \#\{w:\norm{w}\ll U,\ F(w)=0\}
 +\frac1{|\mathscr P(P)|}
 \sum_{\substack{\norm{w}\ll U\\F(w)\ne0}}
 \ \sum_{\substack{\p\in\mathscr P(P)\\p\mid\mathrm N_{\K/\Q}F(w)}}1,
\]
where \(w\) ranges over integral vectors in \(\K^{b+1}\) and \(p\) is the rational prime
below \(\p\).

The first term is \(O(W^{b-1})\) in degree at least five and \(O(W^{b-1/2})\) in
degree four. For a geometrically integral image this follows from the number-field
dimension growth theorem of Paredes and Sasyk \cite[Theorem~1.11]{PS}. More explicitly,
at height at most \(V\) the projective counts are \(O(V^{e(b-1)})\) and
\(O(V^{e(b-1+1/2)})\), respectively, taking \(\varepsilon=1/2\) in degree four.
Balancing coordinates as in \eqref{eq:balanced}, a point of height comparable to
\(V\) has \(O((U/V)^e)\) nonzero integral representatives in the box and requires
\(V\ll U\). The dyadic sums have the asserted size since \(b>2\).
If the image is \(\K\)-irreducible and geometrically reducible, write its reduced
geometric cone as \(\bigcup_\nu Z_\nu\). Galois acts transitively on these distinct
\(b\)-dimensional components. A \(\K\)-point on one \(Z_\nu\) is fixed by Galois
and hence lies on all its conjugates. It belongs to the \(\K\)-defined union of the
pairwise intersections, of dimension at most \(b-1\). Finite linear normalization
on this union bounds the integral \(w\) by \(O(W^{b-1})\). This reasoning concerns
the characteristic-zero condition \(F(w)=0\).

If \(F(w)\ne0\), each prime that contributes divides \(\mathrm N_{\K/\Q}F(w)\), after
a fixed denominator is cleared. This is a nonzero integer of size \(O((1+U)^{e\deg F})\),
with only \(O(1)\) prime divisors in \([P,2P]\) under the hypothesis on \(W\).
There are \(O(W^{b+1})\) projected vectors in all and boundedly many lifts per vector
at each prime. The resulting mean is \(O(W^{b+1}\log P/P)=O(W^{b-1})\).
\end{proof}

\subsection{Completion of the sum}\label{sec:completion}
Let \(N_\K^{\mathrm{ind}}(H)\) count the pairs in Proposition~\ref{pr:independent}. Put \(Q=H^e\).
We leave the size \(P\) of the auxiliary primes free. It is convenient to write
\begin{equation}\label{eq:scales}
 P=QW,\qquad U=W^{1/e},\qquad 2^{e+1}\le W\le Q^{1/2}.
\end{equation}
In this range \(\log P\le\tfrac32\log Q\), so
\(W^2\log P/P=(W/Q)\log P\ll Q^{-1/2}\log Q\). In particular the averaging
hypothesis of Lemma~\ref{le:short} holds.

We work with \(H\) sufficiently large. Let \(B_H\) denote the integral box of
pairs in \((\OO_\K^n)^2\) with norm at most \(H\). It injects modulo each chosen
\(\p\), since a nonzero difference in one coordinate would have absolute norm at
least \(p>(2H)^e\). Take a pair counted by \(N_\K^{\mathrm{ind}}(H)\) and fix one nonzero minor
of its two rows. The minor has absolute norm at most \(2^eQ^2\), so only \(O(1)\)
primes of \(\mathscr P(P)\) can divide it. Write \(\mathcal C^\circ(\F_p)\) for the independent
commuting pairs over \(\F_p\). Counting incidences with primes gives
\begin{equation}\label{eq:sieve}
 (|\mathscr P(P)|-O(1))N_\K^{\mathrm{ind}}(H)\le
 \sum_{\p\in\mathscr P(P)}\sum_{v\in B_H}\ind_{\mathcal C^\circ(\F_p)}(\bar v).
\end{equation}
Since \(|\mathscr P(P)|\gg P/\log P\), division replaces the sum on the right by
\((1+o(1))\) times its prime average.

We use the box weight of \cite[(4.1)]{BSW} in real coordinates on \((\K\otimes_\Q\R)^M\).
Write the \(eM\) coordinates as \(t_1,\ldots,t_{eM}\), taking real and
imaginary parts at the complex places. Let \(\sinc(t)=\sin(\pi t)/(\pi t)\), with
\(\sinc(0)=1\), and put
\[
 w_H(t)=(\pi/2)^{2eM}\prod_{\nu=1}^{eM}\sinc(t_\nu/(2H))^2 .
\]
On \(B_H\) all \(|t_\nu|\le H\), so \(w_H\ge1\). Use Lebesgue measure \(dt\) and
the Fourier transform with negative phase for the trace pairing. If
\(\Tr(\eta\cdot t)=\sum_\nu\eta_\nu^*t_\nu\), the transform is
\[
 \widehat w_H(\eta)=(\pi/2)^{2eM}(2H)^{eM}\prod_{\nu=1}^{eM}(1-2H|\eta_\nu^*|)_+,
 \qquad u_+=\max(u,0).
\]
Indeed the transform of \(\sinc^2\) is the convolution of two unit-interval indicators
centred at the origin. At a complex place the two dual coordinates are
\(2\operatorname{Re}\eta\) and \(-2\operatorname{Im}\eta\).

Choose a generator \(\varpi\) of \(\p\) with \(|\sigma\varpi|\asymp p^{1/e}\) at
every embedding. The dual of \(\p\) for the trace pairing is \(\varpi^{-1}\diff_\K^{-1}\),
where \(\diff_\K\) is the different. Put \(\Delta=\covol(\OO_\K^M)\) and take any lift
\(\widetilde v_0\) of \(v_0\in\mathcal C^\circ(\F_p)\). Summing over its lattice coset
in \(\OO_\K^M\) and then over \(v_0\) gives
\begin{align}
 \sum_{v_0\in\mathcal C^\circ(\F_p)}\sum_{v\equiv v_0\,(\p)}w_H(v)
 &=\frac1{\Delta p^M}\sum_{\xi\in(\diff_\K^{-1})^M}
    \widehat w_H(\xi/\varpi)\,S^\circ_\p(\bar\xi),\label{eq:poisson}\\
 S^\circ_\p(\bar\xi)&=\sum_{v_0\in\mathcal C^\circ(\F_p)}
 \exp\!\big(2\pi i\,\Tr_{\K/\Q}(\xi\cdot\widetilde v_0/\varpi)\big).\notag
\end{align}
For this weight the coset sum converges locally uniformly under translation.
Each unit cube meets a bounded number of points of a fixed lattice, and the
product of the one-dimensional quadratic decays is summable over cubes. Its Fourier
series is finite by the displayed support of \(\widehat w_H\), giving \eqref{eq:poisson}.

The annihilator in the numerator lattice in one coordinate is \(\varpi\diff_\K^{-1}\).
Multiplying the numerator lattice by a fixed integer clearing its denominators
identifies its residue line with \(\F_p\) at every good prime. The resulting additive
character is nontrivial. The inverse Gram matrix of \(\kappa\) converts each of the
two frequency vectors to the coordinates of \eqref{eq:transform}. In coordinate
space we take the pullbacks of \(J,Y,X_j\) by this fixed identification. In particular
it preserves linear dependence. Only \(\norm{\xi}\ll p^{1/e}/H\ll U\) occur in
\eqref{eq:poisson}, and \(H^{eM}/p^M\asymp W^{-M}\). Write
\[
 R_\p=\sum_{\substack{\xi\in(\diff_\K^{-1})^M\\0<\norm{\xi}\le c_1 U}}
 |S^\circ_\p(\bar\xi)|,
\]
where \(c_1\) is a fixed support constant. By \eqref{eq:crude}, \eqref{eq:sieve}
and \eqref{eq:poisson},
\begin{equation}\label{eq:shortsum}
 N_\K^{\mathrm{ind}}(H)\ll W^{-M}\big(P^{n+r}+\Av_\p R_\p\big).
\end{equation}
It remains to estimate this sum of nonzero frequencies.

\subsection{Estimating the short Fourier sum}
Split \(R_\p\) into the contribution of the dependent frequencies, that of the
independent frequencies in \(Y\), and the remainder. Throughout, membership in a
stratum refers to reduction modulo \(\p\). Dependence can be tested over \(\K\).
After fixed denominators are cleared, any nonzero minor of the two frequency vectors
has absolute norm \(O(W^2)=o(P)\) by \eqref{eq:scales}. The same argument in one
coordinate shows that a nonzero frequency in the sum cannot reduce to zero.

\head{Dependent frequencies}
There are \(O(U^{en})=O(W^n)\) choices by the slope count of Section~\ref{sec:heights}.
Except in the directions \(J\), each costs at most \(O(P^{n-1/2})\) by \eqref{eq:dep}.
To count the exceptions, fix a balanced slope of height between \(R\) and \(2R\).
Its parameter vector has norm \(O(U/R)\), and at least one slope coordinate is a unit
at \(\p\) since their joint ideal belongs to a fixed finite set. As \(J\) is a cone,
this vector must itself reduce into \(J\). The first part of Lemma~\ref{le:short}
with \(b=n-1\) therefore bounds the number of bad frequencies by
\[
 \sum_{\substack{R=2^j\ll U\\j\ge0}}R^{2e}(U/R)^{e(n-1)}
 \ll U^{e(n-1)}=W^{n-1},
\]
where \(n\ge8\) and constants absorb the last shell. Thus for each \(\p\)
\[
 R_{\p,\rm dep}\ll W^nP^{n-1/2}+W^{n-1}P^n.
\]

\head{Independent frequencies}
The uniform bound on \(Y\) is only \(O(P^n)\), but Lemma~\ref{le:short} with \(b=n\)
bounds the mean number of frequencies there by \(O(W^{n-1})\). Thus
\(\Av_\p R_{\p,Y}\ll W^{n-1}P^n\). On its complement \eqref{eq:ind} gives \(O(P^{n-1})\),
and \eqref{eq:FK} gives a smaller bound on the shallower strata. In each case
\(S^\circ=S+p\) and the bound absorbs \(p\).

Write \(R_{\p,\rm out}\) for this last sum. On \(X_{k-2}\) apply Lemma~\ref{le:short}
with \(b=M-k+2=n+r+2\). Here \(k-2\ge4\), with equality for \(A_2\). On
\(X_{k-3}\setminus X_{k-2}\) take \(b=n+r+3\), with \(k-3\ge3\), saving at least
\(W^{1/2}\). The spanning hypothesis follows from Lemma~\ref{le:cones}. On
\(X_j\setminus X_{j+1}\) for \(j\le k-4\) the elementary bound \(O(W^{M-j})\)
suffices. Multiplying counts by transform bounds gives
\[
 \begin{aligned}
 \Av_\p R_{\p,\rm out}
 &\ll W^{n+r+1}P^{n-1}+W^{n+r+5/2}P^{n-3/2}
 +\sum_{j=0}^{k-4}W^{M-j}P^{(n+r+j)/2}\\
 &\ll W^{n+r+1}P^{n-1}.
 \end{aligned}
\]
Indeed successive summands in the \(j\)-sum increase by \(\sqrt P/W\).
Its last term is \(W^{n+r+4}P^{n-2}\). Relative to the first term on the right,
the two errors are at most \((W^3/P)^{1/2}\) and \(W^3/P\), respectively.
Both are at most one since \(W^2\le Q\).

The table lists the individual contributions to \eqref{eq:shortsum} after division
by \(Q^n\). Its last four rows count only independent frequencies, with \(Y\) removed
from the last three. All entries in the last three columns are upper bounds up to a
constant factor.
\begin{center}
\small\setlength{\tabcolsep}{4pt}
\begin{tabular}{@{}lccc@{}}
\toprule
Frequency & Mean count & \(\lvert S^\circ_\p\rvert\) & Relative contribution\\
\midrule
Zero & \(1\)& \(P^{n+r}\)&\(Q^r/W^k\)\\
Dependent, direction outside \(J\)& \(W^n\)&\(P^{n-1/2}\)&\(P^{-1/2}\)\\
Dependent, direction in \(J\)& \(W^{n-1}\)&\(P^n\)&\(W^{-1}\)\\
\(Y\)&\(W^{n-1}\)&\(P^n\)&\(W^{-1}\)\\
\(X_{k-2}\setminus Y\)&\(W^{n+r+1}\)&\(P^{n-1}\)&\(W^{r+1}/P\)\\
\(X_{k-3}\setminus X_{k-2}\)&\(W^{n+r+5/2}\)&\(P^{n-3/2}\)&\(W^{r+5/2}/P^{3/2}\)\\
\(X_j\setminus X_{j+1}\), \(0\le j\le k-4\)&\(W^{M-j}\)&\(P^{(\dcom+j)/2}\)&\(W^{n-j}/P^{(k-j)/2}\)\\
\bottomrule
\end{tabular}
\end{center}
\head{Choice of the prime scale}
Substitute the three sums into \eqref{eq:shortsum} and use \(P^{-1/2}\le W^{-1}\).
With \(P=QW\) and \(M=2n\) this gives
\begin{equation}\label{eq:optimization}
 \begin{aligned}
 N_\K^{\mathrm{ind}}(H)&\ll W^{-2n}\big(P^{n+r}+W^{n-1}P^n+W^{n+r+1}P^{n-1}\big)\\
 &\ll \frac{Q^{n+r}}{W^{n-r}}+\frac{Q^n}{W}+Q^{n-1}W^r .
 \end{aligned}
\end{equation}
The first two terms decrease with \(W\) while the third increases. The third
overtakes the first at \(W^n=Q^{r+1}\) and the second at \(W^{r+1}=Q\). We therefore take
\[
 W=\max\{Q^{(r+1)/n},Q^{1/(r+1)}\}.
\]
Since \(r\ge2\) and \(n\ge r(r+2)\), this choice is at most \(Q^{3/8}\) and satisfies
\eqref{eq:scales} for large \(H\). Now \eqref{eq:optimization} gives
\[
 N_\K^{\mathrm{ind}}(H)\ll Q^{n-1}W^r=Q^{n-\min\{1-r(r+1)/n,\,1/(r+1)\}}=H^{en-e\delta(n,r)}.
\]
Increasing the implicit constant covers bounded \(H\). This proves
Proposition~\ref{pr:independent}, subject only to Proposition~\ref{pr:geometry}.

For example, in \(\slie_3(\Q)\) one has \(Q=H,\ n=8\) and \(r=2\). The estimate
before choosing the prime is
\[
 N_\Q^{\mathrm{ind}}(H)\ll H^{10}W^{-6}+H^8W^{-1}+H^7W^2.
\]
The first and last terms agree at \(W=H^{3/8}\), so \(P=H^{11/8}\).
They are both \(H^{31/4}\), while the middle term is \(H^{61/8}\).
Thus independent commuting pairs contribute \(O(H^{31/4})\). The proportional pairs
supply the main term \(\frac{3\zeta(7)}{\zeta(8)}(2\lfloor H\rfloor+1)^8\) for the trace-zero entry box by
Corollary~\ref{cor:simple} and \eqref{eq:tau}.

\section{Geometry of the bracket}\label{sec:geom}
We prove Proposition~\ref{pr:geometry}. In this section \(\s\) is simple over an
algebraically closed field of characteristic zero, \(G\) is its connected adjoint
group, and \(\s^y\) denotes the centralizer of \(y\). Recall that \(F_a\) is the
fibre \([y,u]=a\), and that \(F_{a,b}\) is its intersection with \(\kappa(b,y)=0\).

For fixed \(y\), invariance of \(\kappa\) gives \([y,\s]=(\s^y)^\perp\). Hence
\[
 \{u\in\s:[y,u]=a\}=
 \begin{cases}
  u_0+\s^y,& a\perp\s^y,\ \ [y,u_0]=a,\\
  \varnothing,& a\not\perp\s^y.
 \end{cases}
\]
On each Jordan class \(\dim\s^y\) is constant. It remains to bound the dimension
of the locus \(a\perp\s^y\) in that class.

\subsection{Projection of an orbit onto a torus}
For a torus \(T\subset G\) write \(\tor=\Lie T\), \(t=\dim T\) and
\[
 \mathcal I_T(a)=\{g\in G:\Ad(g^{-1})a\perp\tor\}.
\]
All orthogonals here are taken with respect to \(\kappa\).

\begin{lemma}\label{le:torus}
For \(a\ne0\), the part of \(\mathcal I_T(a)\) on which \(\dim\tor^{\Ad(g^{-1})a}=j\)
has dimension at most \(n-t-j\).
\end{lemma}
\begin{proof}
On the orbit \(Ga\) stratify \(v\) by the support of its \(T\)-weight decomposition.
Then \(\mathfrak u=\tor^v\) is fixed on each stratum. Put
\(\mathfrak m=\s^{\mathfrak u}\) and \(T_v(Ga)=[\s,v]\). The representation of
\(\mathfrak u\) on \(T_v(Ga)\) is faithful. Indeed, if \(0\ne h\in\mathfrak u\)
kills \(T_v(Ga)\), then \(v\) commutes with each nonzero eigenspace of \(\ad h\),
since \([v,x]\) for an eigenvector \(x\) there has the same weight and belongs to
\(T_v(Ga)\). The subalgebra generated by those eigenspaces is stable under the zero
eigenspace and contains every nonzero eigenspace. It is a nonzero ideal, and equals
\(\s\). This would make \(v\) central.

The orbit form \(\omega([x,v],[y,v])=\kappa(v,[x,y])\) is alternating, nondegenerate
and \(\mathfrak u\)-invariant. Faithfulness and the pairing of opposite weights show
that at least \(2j\) dimensions of \(T_v(Ga)\) have nonzero weight. Its fixed
subspace is \(T_v(Ga)^{\mathfrak u}=[\mathfrak m,v]\). Projection from this space to
\(\tor^*\) has rank \(t-j\). In fact an element \(h\) of \(\tor\) kills its image
precisely if \(\kappa([v,h],\mathfrak m)=0\), and \([v,h]\in\mathfrak m\) with
\(\kappa\) nondegenerate on \(\mathfrak m\). Every point of the support stratum is
fixed by \(\mathfrak u\), so the stratum lies in \(\mathfrak m=\s^{\mathfrak u}\).
At a smooth point its tangent space lies in
\(T_v(Ga)\cap\mathfrak m=T_v(Ga)^{\mathfrak u}\). Cutting also by \(\tor^\perp\)
places the tangent space of the intersection in the kernel of the projection.
Adding the fibre dimension of the orbit map, the required upper bound is
\[
 (\dim Ga-2j)-(t-j)+(n-\dim Ga)=n-t-j.\qedhere
\]
\end{proof}
\begin{proof}[Proof of the fibre dimension assertions]
Use the finite stratification into Jordan classes. In one such class write
\(y=\Ad(g)(y_s+v)\) with \(\lev=\s^{y_s}\) a fixed Levi subalgebra, \(v\in\lev\)
a chosen representative of a fixed nilpotent orbit, and \(y_s\) in the open stratum
of \(\mathfrak c=Z(\lev)\) with centralizer \(\lev\). Let \(t=\dim\mathfrak c\),
\(c_v=\dim\lev^v\), and \(L_0\) the connected Levi group.
The condition \(a\in[y,\s]\) is equivalent to \(\Ad(g^{-1})a\perp\lev^v\).
This defines a closed subset \(A\) of \(\mathcal I_{Z(L_0)^\circ}(a)\) of dimension
at most \(n-t\), for \(a\ne0\). For each admissible \(y\), the \(u\)-fibre has dimension
\(c_v\). The element centralizer in \(L_0\) contributes at least \(c_v\)
dimensions of redundancy among the parameters \((g,y_s)\), since it fixes \(y_s+v\).
If \(F_a^{\rm cl}\) denotes the part of the fibre above this Jordan class, we obtain
\[
 \dim F_a^{\rm cl}\le(\dim A+t-c_v)+c_v\le(n-t)+t=n.
\]

If the contribution reaches \(n\), then \(A\) contains an irreducible component
\(I_0\) of \(\mathcal I_{Z(L_0)^\circ}(a)\) of dimension \(n-t\).
Right multiplication by \(L_0\) preserves \(I_0\) by connectedness. Consequently
\(\Ad(g^{-1})a\) is perpendicular to the span of \(\Ad(L_0)\lev^v\) for every
\(g\in I_0\). The span contains the centre of \(\lev\) and a nonzero vector in each
simple ideal, by taking the centralizer of the corresponding component of \(v\).
Its \(L_0\)-invariance then gives the whole of \(\lev\). Hence \(I_0\) lies in
\(\mathcal I_T(a)\) for a maximal torus \(T\) of \(L_0\), of rank \(r\).
Lemma~\ref{le:torus} gives \(n-t\le n-r\), so \(t=r\) and \(y\) is regular
semisimple. Every component of a nonempty fibre has dimension at least \(n\) by
the height theorem. This proves purity and the density assertion. When \(a=0\)
there is no condition on \(g\) and the same calculation gives \(n+t\), with
equality \(n+r\) on regular semisimple elements. Richardson's theorem also gives
irreducibility at zero~\cite{Richardson}.

Surjectivity is classical~\cite{Brown}. It also follows here from the dimension
calculation. Miracle flatness gives flatness over the punctured target, so the image
there is open. A nonempty complement would project to a closed \(G\)-stable subset
of \(\PP(\s)\), by scalar invariance. The Borel fixed point theorem would put the
highest-root line in this subset, whereas \(e_\theta=[h,e_\theta]\) for
\(\theta(h)=1\).
\end{proof}
For a reductive algebra over the same field, write the geometric simple ideals as
\(\s'_j\) with ranks \(r'_j\). Factoring the bracket into the product of the simple
bracket maps gives, for \(a=(a_j)\) in the derived algebra,
\begin{equation}\label{eq:reductivefibre}
 \dim\{(x,y):[x,y]=a\}=d+z+\sum_{j:a_j=0}r'_j .
\end{equation}
In particular the largest open subset of the derived algebra over which the bracket
is flat consists of points whose every simple component is nonzero. This includes
the abelian case with a one-point target.

\subsection{Irreducibility for a general target}
\begin{proof}[Proof of Proposition~\ref{pr:geometry}(i)]
Let \(f:X=\s^2\to\s\) be the bracket. The equations of \(X\times_\s X\) are
\([x_1,y_1]=[x_2,y_2]\). Replacing \(y_2\) by \(-y_2\) identifies this variety
with the zero fibre of the sum of two brackets. Glazer and Hendel
\cite[Ch.~6, pp.~67--70]{GH} prove that this sum is flat with geometrically
irreducible, reduced fibres. In particular \(X\times_\s X\) is irreducible. Its
map to \(\s\) is dominant, so over the function field its generic fibre
\(X_\eta\times_{k(\s)}X_\eta\) is irreducible as well.

If the geometric generic fibre of \(f\) had components \(C_1,\ldots,C_q\) with
\(q>1\), the components of its square would be \(C_i\times C_j\). The union over
\(i=j\) and the union over \(i\ne j\) would be two proper closed subsets, each
stable under the Galois action. They descend to the function field and cover the
generic square, a contradiction. Geometric irreducibility thus holds over a dense
open subset by constructibility~\cite{Stacks}. The closure of the failure locus,
with zero adjoined, is a proper cone by homogeneity.
\end{proof}

\subsection{Hyperplanes containing a component}
The exceptional set in Proposition~\ref{pr:geometry}(ii) can be nonempty. In type
\(A_2\), take \(a=E_{12}\) and \(b=E_{32}\) in \(\slie_3\). On the six-dimensional
linear space \(y_{21}=y_{23}=0\), regular semisimple elements form a nonempty open
subset (it contains regular diagonal matrices). Here \((y^2)_{21}=0\). The centralizer
of such a \(y\) is spanned by \(y\) and the traceless part of \(y^2\), so \(a\)
annihilates it under the trace form. Thus \([y,u]=a\) has an affine plane of solutions.
Moreover \(\Tr(by)=y_{23}=0\). Hence \(\dim F_{a,b}=8=\dim\slie_3\), although
\(a,b\) are independent.

\begin{proof}[Proof of Proposition~\ref{pr:geometry}(ii)]
Fix a Cartan \(\h\) with torus \(T\), and put \(\mathcal I_a=\mathcal I_T(a)\).
Suppose \(\dim F_{a,b}=n\) with \(a,b\) independent. An irreducible component \(C\)
of \(F_{a,b}\) of dimension \(n\) is closed in \(F_a\). As every component of
\(F_a\) has dimension \(n\), \(C\) is a component of \(F_a\). Hence the locus
\(C^{\rm rs}\) where \(y\) is regular semisimple is dense in \(C\) by the first
part of the proposition.

Map \((y,u)\in C^{\rm rs}\) to the unique Cartan \(\s^y\) containing \(y\).
Since \(a=[y,u]\), this Cartan is perpendicular to \(a\). These Cartans are
parametrized by \(\mathcal I_a/N_G(T)\), of dimension at most \(n-2r\) by
Lemma~\ref{le:torus}. Over a fixed Cartan \(\h'\) the choice of \(y\) has
dimension at most \(r\), and the solutions \(u\) form an affine space under \(\h'\)
of dimension \(r\). Thus the closure of the image of \(C^{\rm rs}\) has dimension
exactly \(n-2r\), and its general fibre has dimension \(2r\). If \(\h'\) is not
perpendicular to \(b\), the condition \(\kappa(b,y)=0\) cuts the allowed \(y\) to
dimension at most \(r-1\). The general Cartan in the image is therefore perpendicular
to \(b\). Orthogonality to \(b\) is closed on the variety of Cartans, so this holds
on the entire image closure there. Its inverse image in \(\mathcal I_a\) has dimension
\((n-2r)+r=n-r\), because the normalizer has dimension \(r\). It lies in \(\mathcal I_b\)
and contains a component of \(\mathcal I_a\) of dimension \(n-r\).

Write \(a'=\Ad(g^{-1})a,\ b'=\Ad(g^{-1})b\). Lemma~\ref{le:torus} puts the locus
\(\dim\h^{a'}>0\) in dimension at most \(n-r-1\). On a dense open subset of the
component one therefore has \(\h^{a'}=0\). There the \(r\) equations
\(\kappa(a',h)=0\) (\(h\) in a basis of \(\h\)) have independent differentials,
corresponding under \(\kappa\) to \([a',h]\). Their common kernel is the tangent
space of the smooth component. Since \(\kappa(b',h)=0\) on this component, its
differentials vanish on the same tangent space. We have proved
\[
 a',b'\in\h^\perp,\qquad \h^{a'}=0,\qquad [b',\h]\subseteq[a',\h]
\]
for a set of \(g\) of dimension \(n-r\) above each exceptional pair.

Let \(S\) be the support of \(a'\) in the root system \(\Phi\). It spans \(\h^*\)
since \(\h^{a'}=0\). There is an endomorphism \(U\) of \(\h\) with
\([b',h]=[a',Uh]\) for every \(h\in\h\). Comparing root coordinates shows that
\(b'\) is supported on \(S\), and each \(\alpha\in S\) is an eigenvector of \(U^*\)
with eigenvalue \(b'_\alpha/a'_\alpha\). These roots span, so \(U^*\) is
diagonalizable. For an independent pair there are \(q\ge2\) eigenvalues.
Their eigenspaces have underlying real root spans \(V_1,\ldots,V_q\), whose
dimensions sum to \(r\).

Let \(h_\Phi=|\Phi|/r\), the Coxeter number. In the real reflection representation
the sum of the orthogonal projections onto the root lines, counting both signs of
a root, is \(h_\Phi\) times the identity by Weyl invariance and trace. For a nonzero
proper root span \(V\) it follows that
\[
 |\Phi\cap V|\le h_\Phi\dim V-1.
\]
Indeed equality with \(h_\Phi\dim V\) would make every root outside \(V\)
orthogonal to \(V\), contradicting irreducibility. The strict inequality saves at
least one by integrality. Summing over \(V_j\) gives \(|S|+q\le|\Phi|=n-r\).
For fixed \(g\), support and partition of the support, at most \(|S|+q\) parameters
determine the root coordinates and eigenvalues. There are finitely many supports
and partitions. To pass from this count to the exceptional set, let \(\Xi\)
be the set of triples \((g,a,b)\in G\times\s^2\) for which \(a,b\) are independent,
\(a',b'\perp\h,\ \h^{a'}=0\) and \([b',\h]\subseteq[a',\h]\). These are
constructible rank conditions on the two maps \(\h\to\s\). For a fixed \(g\)
the preceding count gives \(\dim\Xi_g\le|\Phi|=n-r\), so \(\dim\Xi\le2n-r\).
Above each exceptional pair the dense open subset just constructed lies in the
fibre of \(\Xi\to\s^2\). This fibre has dimension at least \(n-r\). The
exceptional set is constructible by the fibre-dimension theorem, and its closure
therefore has dimension at most
\[
 \dim\Xi-(n-r)\le(2n-r)-(n-r)=n,
\]
as asserted.
\end{proof}

\section{Remarks on the error term}\label{sec:remainder}
Two different estimates determine the remainder in \eqref{eq:error}. The first
concerns the boundary of the slope sections. Even in dimension one a term of this
size can depend on the dilation. For \(\g=\Q,\ \Lambda=\Z\) and \(B=[-1,1]\),
every pair commutes and
\[
 N_{\Z,[-1,1]}(H)-4H^2=4(1-2\{H\})H+O(1).
\]
Along integer \(H\) the coefficient of \(H\) is \(4\), whereas along half-integer
\(H\) it is zero. Thus the boundary estimate alone cannot supply a constant second
coefficient for real \(H\). In the \(L\)-normalization this example gives \(N=L^2\)
exactly.

Take \(\Lambda=\slie_2(\Z)\) and let \(B\) be the unit cube in a lattice basis.
Then \(E=3\). Write \(N(h)=N_{\Lambda,B}(h)\) for integers \(h\ge0\).
Each nonzero commuting pair has a unique slope. Grouping the \(4\varphi(j)\) slopes
of height \(j\) gives
\[
 N(h)=1+4\sum_{1\le j\le h}\varphi(j)
       \bigl((2\lfloor h/j\rfloor+1)^3-1\bigr).
\]
Only \(j\mid h\) contribute to the increment at \(h\ge1\). Since
\((2k+1)^3-(2k-1)^3=24k^2+2\) and \(\sum_{j\mid h}\varphi(j)=h\),
\[
 N(h)-N(h-1)=96\sum_{j\mid h}\varphi(j)(h/j)^2+8h.
\]
For \(R(h)=N(h)-\frac{4\zeta(2)}{\zeta(3)}(2h+1)^3\), this gives
\[
 R(h)-R(h-1)
 =96h^2\left(\sum_{j\mid h}\frac{\varphi(j)}{j^2}-\frac{\zeta(2)}{\zeta(3)}\right)+8\left(h-\frac{\zeta(2)}{\zeta(3)}\right).
\]
At a prime \(q\) the divisor sum is \(1+(q-1)/q^2\). Since \(\zeta(2)/\zeta(3)>1\),
the increment of \(R\) is \(96\bigl(1-\zeta(2)/\zeta(3)\bigr)q^2+O(q)\).

Larger increments occur at \(h=\operatorname{lcm}(1,\ldots,k)\). Every
\(j\le k\) divides \(h\). Writing \(\mu\) for the Möbius function and using
\(\varphi(j)/j=\sum_{d\mid j}\mu(d)/d\), we obtain
\[
 \sum_{j\le k}\frac{\varphi(j)}{j^2}
 =\sum_{d\le k}\frac{\mu(d)}{d^2}\sum_{m\le k/d}\frac1m
 =\frac{\log k}{\zeta(2)}+O(1).
\]
The last equality follows from the harmonic-sum estimate, the absolute
convergence of \(\sum_d(\log d)/d^2\), and
\(\sum_d\mu(d)/d^2=1/\zeta(2)\). Since \(h\le k!\le k^k\),
we have \(\log\log h\ll\log k\). For all sufficiently large \(k\) the
increment identity therefore gives
\[
 \max\{|R(h)|,\ |R(h-1)|\}
 \ge \tfrac12|R(h)-R(h-1)|\gg h^2\log\log h.
\]
Thus \(|R(H)|\gg H^2\log\log H\) along a sequence of integer heights for
the unit lattice-basis cube in \(\slie_2\). In particular, a fixed multiple
of \(L^2\) cannot give a second term with remainder \(o(L^2)\). Whether the
order \(h^2\log\log h\) is a uniform upper bound remains open. For this cube
the following estimate improves the general bound in \eqref{eq:simple}.

\begin{proposition}\label{pr:cube-error}
For the unit lattice-basis cube in \(\slie_2(\Z)\), the remainder at integer
heights \(h\ge1\) satisfies
\[
 R(h)=O\!\left(h^2(\log(2h))^{2/3}\right).
\]
The bound \(R(h)=O(h^2\log\log h)\) as \(h\to\infty\) is equivalent to
\[
 \sum_{m\le x}\sigma(m)=\frac{\zeta(2)}{2}x^2+O(x\log\log x)
 \qquad(x\to\infty),
\]
where \(\sigma(m)=\sum_{d\mid m}d\).
\end{proposition}
\begin{proof}
Substitute \(\varphi(j)=\sum_{d\mid j}\mu(d)j/d\) in the increment of \(N\).
Writing \(j=dk\) and using \(\sum_{k\mid m}1/k=\sigma(m)/m\) gives
\[
 \sum_{j\mid h}\varphi(j)(h/j)^2
 =\sum_{d\mid h}\mu(d)(h/d)\sigma(h/d).
\]
Put \(T(x)=\sum_{m\le x}m\sigma(m)\) and \(A(x)=\sum_{d\le x}\mu(d)T(x/d)\).
Summing the increments from \(1\) to \(h\), including the term \(8\sum_{j\mid n}\varphi(j)=8n\)
at each height \(n\), gives
\[
 N(h)=1+4h(h+1)+96A(h),
\]
since \(N(0)=1\). Subtracting \(\frac{4\zeta(2)}{\zeta(3)}(2h+1)^3\) yields
\(R(h)=96\bigl(A(h)-\tfrac{\zeta(2)}{3\zeta(3)}h^3\bigr)+O(h^2)\).

Write \(S(x)=\sum_{m\le x}\sigma(m)\), and let \(g(x)\) be nondecreasing and
bounded below by a positive constant on \([1,\infty)\). Partial summation gives
\[
 T(x)=xS(x)-\int_1^x S(t)\,dt,\qquad
 S(x)=\frac{T(x)}x+\int_1^x\frac{T(t)}{t^2}\,dt .
\]
Thus \(S(x)=\zeta(2)x^2/2+O(xg(x))\) is equivalent to
\(T(x)=\zeta(2)x^3/3+O(x^2g(x))\).
Inserting the latter estimate into \(A(x)\) gives
\(A(x)=\tfrac{\zeta(2)}{3\zeta(3)}x^3+O(x^2g(x))\), since
\(\sum_{d\ge1}\mu(d)/d^3=1/\zeta(3)\) and \(\sum_{d\ge1}d^{-2}<\infty\).
Conversely, Möbius inversion gives \(T(x)=\sum_{d\le x}A(x/d)\), recovering
the estimate for \(T\) by the same summation. In both directions extending the
cubic main sum to all \(d\) costs \(O(x)\). Since \(A(x)=A(\lfloor x\rfloor)\)
and \(x^3-\lfloor x\rfloor^3=O(x^2)\), the estimate for \(A\) is equivalent to
\(R(h)=O(h^2g(h))\) at integer heights.

Walfisz's estimate for the sum-of-divisors function \cite[Satz~4, p.~99]{Walfisz},
also stated by Tóth \cite[proof of Theorem~22, p.~17]{Toth}, is
\[
 \sum_{m\le x}\sigma(m)=\frac{\pi^2}{12}x^2+O\!\left(x(\log x)^{2/3}\right)
 \qquad(x\to\infty).
\]
Since \(\pi^2/12=\zeta(2)/2\), taking \(g(x)=(\log(2x))^{2/3}\) proves the
cube bound. Taking \(g(x)=\log\log(3x)\) gives the stated equivalence.
\end{proof}

For the divisor-sum error \(D(x)=S(x)-\zeta(2)x^2/2\), Pétermann's theorem
and addendum~\cite{Petermann} give the two-sided estimate
\(D(x)=\Omega_{\pm}(x\log\log x)\) (see \cite[pp.~171--172]{PetermannFluctuations}
for this normalization). By Proposition~\ref{pr:cube-error}, the open upper-bound
question for the cube is precisely whether \(\lvert D(x)\rvert\ll x\log\log x\)
for all sufficiently large \(x\).

The other estimate concerns independent pairs. In the notation of
\eqref{eq:optimization}, the cost of the deepest Fourier strata is \(Q^{n-1}W^r\).
In type \(A_r\) one has \(n=r(r+2)<(r+1)^2\). At the optimal scale this term has
the same size as the zero-frequency term \(Q^{n+r}/W^{n-r}\). In the remaining
types of rank at least two, \(n\ge(r+1)^2\) and the balancing term is \(Q^n/W\),
from \(Y\) and the exceptional dependent directions. This accounts for the two
different savings in \eqref{eq:delta}. A rank-one simple factor contributes only
to the lattice estimate, since all its commuting pairs are dependent.

\appendix
\section{Box heights for simple and similitude Lie algebras}\label{sec:boxes}
All entries in the tables below are coefficients of \(L^E\), where
\(L=2\lfloor H\rfloor+1\) and \(E=\dim\g+\dim Z(\g)\). We first bound the coordinates in a lattice basis.
We then fix a matrix representation and bound every matrix entry. The distinction
already occurs in \(\slie_3\). Bounding two independent diagonal entries gives a
square, while bounding all three removes two corners of that square.

\subsection{The unit box in a lattice basis}
Choose a full lattice \(\Lambda\) and a \(\Z\)-basis \(\mathcal B=(v_1,\ldots,v_n)\).
For \(x=\sum_i u_i v_i\), define \(h_{\mathcal B}(x)=\max_i|u_i|\) and
\(B_{\mathcal B}=\{x:h_{\mathcal B}(x)\le1\}\). The map \(u\mapsto\sum_i u_i v_i\)
identifies \(\Z^n\) with \(\Lambda\) and \([-1,1]^n\) with \(B_{\mathcal B}\).
Thus \(\vol_\Lambda(B_{\mathcal B})=2^n\). In absolute rank one commutation is
equivalent to linear dependence, so the basis identifies the count with the cube
count of Proposition~\ref{pr:cube-error} at \(h=\lfloor H\rfloor\).
Together with Corollary~\ref{cor:simple} in higher rank, this gives for an
absolutely simple algebra
\begin{equation}\label{eq:basis}
 N_{\Lambda,B_{\mathcal B}}(H)=\frac{4\zeta(n-1)}{\zeta(n)}L^n+
 \begin{cases}O\!\left(H^2(\log(2H))^{2/3}\right),&r=1,\\
 O_{\Lambda,\mathcal B}(H^{n-\delta(n,r)}),&r\ge2.\end{cases}
\end{equation}
For a split simple algebra this determines the lattice-basis coefficient in every
type. In the classical types \(A_r,B_r,C_r,D_r\), substitute
\(n=r(r+2),r(2r+1),r(2r+1),r(2r-1)\), respectively. In types \(G_2,F_4,E_6,E_7,E_8\)
the respective coefficients are
\[
 \frac{4\zeta(13)}{\zeta(14)},\quad \frac{4\zeta(51)}{\zeta(52)},\quad
 \frac{4\zeta(77)}{\zeta(78)},\quad \frac{4\zeta(132)}{\zeta(133)},\quad
 \frac{4\zeta(247)}{\zeta(248)}.
\]

\head{Similitudes}
For a nonsingular alternating or symmetric \(D\times D\) matrix \(J\), take the
full integral lattice
\[
 \Lambda_J=\{X\in M_D(\Z):X^{\mathsf T}J+JX=\lambda J\text{ for some }\lambda\in\Q\}.
\]
We write \(\gsp_{2r}\) for the standard alternating form and \(\gso_D\) for the split
symmetric form \(J_D\) with ones on the antidiagonal. When the derived algebra is
absolutely simple over \(\Q\) and has dimension \(m\), choose a lattice basis starting
with \(I_D\). This is possible since \(\Z I_D\) is primitive. The remaining basis vectors project to a
basis of \(\bar\Lambda=\Lambda_J/\Z I_D\), with coordinate cube \(\bar B\).
There are exactly \(L\) scalar lifts of each quotient point in \(H\bar B\) to the
dilated coordinate box. The bracket lies in the derived algebra, which meets the
scalar line trivially, so commutation is detected in the quotient. Consequently
\(N_{\Lambda_J,B_{\mathcal B}}(H)=L^2N_{\bar\Lambda,\bar B}(H)\). By \eqref{eq:basis}
the coefficient of \(L^{m+2}\) is \(4\zeta(m-1)/\zeta(m)\). The two-factor case in split dimension
four is given below.

\subsection{The entrywise box in a matrix representation}
Fix a faithful representation \(\rho:\g\to M_D(\Q)\) and a basis of its representation
space. Bound the entries of the representing matrix by taking
\[
 \Lambda_\rho=\{x\in\g:\rho(x)\in M_D(\Z)\},\qquad
 B_\rho=\{x\in\g_\R:\max_{i,j}|\rho(x)_{ij}|\le1\}.
\]
When \(\g\) is absolutely simple over \(\Q\), Corollary~\ref{cor:simple} gives,
for any lattice basis \(\mathcal B\) of \(\Lambda_\rho\), the coefficient of \(L^n\) as
\begin{equation}\label{eq:representation}
 K_\rho=\frac{4\zeta(n-1)}{\zeta(n)}\frac{\vol_{\Lambda_\rho}(B_\rho)}{\vol_{\Lambda_\rho}(B_{\mathcal B})}
 =\frac{4\zeta(n-1)}{\zeta(n)}\,2^{-n}\vol_{\Lambda_\rho}(B_\rho).
\end{equation}
The arithmetic factor \(4\zeta(n-1)/\zeta(n)\) is thus common to both heights. It remains to find
the ratio of the two volumes in the same lattice coordinates. Formula \eqref{eq:simple}
gives the error term for either body.

\subsubsection{The defining representations of the classical types}
Use the standard alternating form for \(\spLie_{2r}\) and the antidiagonal matrix
\(J_D\) with entries one on the antidiagonal for the orthogonal algebra in dimension
\(D\). Put
\begin{equation}\label{eq:tau}
 \tau_N=\frac1{(N-1)!}\sum_{k=0}^{\lfloor N/2\rfloor}(-1)^k
 \binom Nk\left(\frac N2-k\right)^{N-1}\qquad(N\ge2).
\end{equation}
Then the four defining representations give the following coefficients.
\begin{center}
\begin{tabular}{@{}llcc@{}}
\toprule
Type & Integral matrix model & \(n\) & \(K_\rho\)\\
\midrule
\(A_{N-1},\ N\ge2\)&\(\slie_N(\Z)\)& \(N^2-1\)&\(4\tau_N\,\zeta(n-1)/\zeta(n)\)\\
\(B_r,\ r\ge2\)&\(\soLie(J_{2r+1})\cap M_{2r+1}(\Z)\)&\(r(2r+1)\)&\(4\zeta(n-1)/\zeta(n)\)\\
\(C_r,\ r\ge2\)&\(\spLie_{2r}(\Z)\)&\(r(2r+1)\)&\(4\zeta(n-1)/\zeta(n)\)\\
\(D_r,\ r\ge4\)&\(\soLie(J_{2r})\cap M_{2r}(\Z)\)&\(r(2r-1)\)&\(4\zeta(n-1)/\zeta(n)\)\\
\bottomrule
\end{tabular}
\end{center}
\head{Trace zero}
In \(\slie_N\) take the \(N(N-1)\) off-diagonal entries and the first \(N-1\)
diagonal entries as free integral coordinates. The last diagonal entry is minus the
sum of the others. After dividing the total volume by \(2^{N^2-1}\), the off-diagonal
factors disappear, leaving
\[
 2^{-(N^2-1)}\vol_{\Lambda_\rho}(B_\rho)=\Pr\{|U_1+\cdots+U_{N-1}|\le1\},
\]
where the \(U_i\) are independent uniform variables on \([-1,1]\). Writing
\(U_i=2T_i-1\), this is the probability that \(\sum T_i\) lies between \(N/2-1\)
and \(N/2\). If \(m=N-1\), inclusion--exclusion in the unit cube gives the distribution
function \(\sum_{k=0}^m(-1)^k\binom mk(s-k)_+^m/m!\) at \(s\).
Subtracting its values at the two endpoints and using Pascal's identity gives
\eqref{eq:tau}. In particular \(\tau_2=1,\ \tau_3=3/4\) and \(\tau_4=2/3\).
For \(\slie_3(\Z)\) the two free diagonal coordinates \(u,v\) range over \([-1,1]^2\)
in the lattice-basis box. The last entry imposes \(|u+v|\le1\) in the entrywise
box, leaving three quarters of the square. The respective main terms are therefore
\(\frac{4\zeta(7)}{\zeta(8)} L^8\) and \(\frac{3\zeta(7)}{\zeta(8)} L^8\), both with error \(O(H^{31/4})\).

\head{Alternating and symmetric forms}
In the standard symplectic basis one has
\[
 x=\begin{pmatrix}D_0&S\\T&-D_0^{\mathsf T}\end{pmatrix},
 \qquad S=S^{\mathsf T},\quad T=T^{\mathsf T}.
\]
All entries of \(D_0\) and the two upper triangles of \(S,T\) are independent integral
coordinates. There are \(r^2+r(r+1)=r(2r+1)\). The matrix entry bound is exactly
the unit cube in these coordinates. For \(x^{\mathsf T}J_D+J_Dx=0\) one has
\(x_{ij}=-x_{D+1-j,D+1-i}\). The \(D\) antidiagonal entries vanish, and one
coordinate from each of the \((D^2-D)/2\) remaining pairs is free and integral.
Its partner has the same absolute value. Both calculations give
\(2^{-n}\vol_{\Lambda_\rho}(B_\rho)=1\). The latter calculation also applies in
dimensions three and six, and the symplectic calculation applies in dimension two.
For example \(D_3=A_3\), but the \(6\)-dimensional orthogonal realization gives
\(\frac{4\zeta(14)}{\zeta(15)}\) and the \(4\)-dimensional trace-zero realization gives
\(\frac{8\zeta(14)}{3\zeta(15)}\). This illustrates the effect of specifying the matrix height.

\subsubsection{Symplectic and split orthogonal similitudes}\label{sec:similitudes}
For the full integral lattices \(\Lambda_J\), the table gives the coefficient of
\(L^E\) for each height. Here \(m\) is the dimension of the derived algebra and
\(E=m+2\). The basis starts with \(I_D\) and, for \(\gso_4\), uses the two quotient
triples specified below.
\begin{center}
\renewcommand{\arraystretch}{1.4}
\begin{tabular}{@{}lccc@{}}
\toprule
Lie algebra & \(m\) & Lattice-basis box & Entrywise box\\
\midrule
\(\gsp_{2r},\ r\ge1\)&\(r(2r+1)\)&\(\dfrac{4\zeta(m-1)}{\zeta(m)}\)&\(\dfrac{4(r+4)\zeta(m-1)}{(r+1)(r+2)\zeta(m)}\)\\\addlinespace[4pt]
\(\gso_{2r+1},\ r\ge1\)&\(r(2r+1)\)&\(\dfrac{4\zeta(m-1)}{\zeta(m)}\)&\(\dfrac{2(r+4)\zeta(m-1)}{(r+1)(r+2)\zeta(m)}\)\\\addlinespace[4pt]
\(\gso_{2r},\ r\ge3\)&\(r(2r-1)\)&\(\dfrac{4\zeta(m-1)}{\zeta(m)}\)&\(\dfrac{4(r+4)\zeta(m-1)}{(r+1)(r+2)\zeta(m)}\)\\\addlinespace[4pt]
\(\gso_4\)&6&\(\bigl(4\zeta(2)/\zeta(3)\bigr)^2\)&\(K_4\) in \eqref{eq:gso4}\\
\bottomrule
\end{tabular}
\end{center}
For the first three rows with \(r\ge2\), \eqref{eq:error} gives
\(O(H^{E-\delta(m,r)})\) for both boxes. In rank one (\(\gsp_2\) and \(\gso_3\)),
the same estimate gives an entrywise error of \(O(H^4\log(2H))\). For the
lattice-basis box, each commuting pair in the quotient has exactly \(L^2\) scalar
lifts, so the rank-one error in \eqref{eq:basis} gives
\(O(H^4(\log(2H))^{2/3})\). The entrywise coefficient for \(\gsp_2=\glie_2\) is
\(\frac{10\zeta(2)}{3\zeta(3)}\), as in Corollary~\ref{cor:matrix}.

Put \(\gamma=\lambda/2\). The diagonal entries occur in pairs
\(\gamma+u_i,\gamma-u_i\) for \(1\le i\le r\), with an additional entry \(\gamma\)
in odd dimension. In the symplectic model the diagonal blocks are \(D_0\) and
\(2\gamma I_r-D_0^{\mathsf T}\), and the other two blocks are symmetric.
In the orthogonal model, reflection across the antidiagonal pairs off-diagonal
entries with opposite signs and forces the fixed ones to vanish. In all three cases
the \(m-r\) free off-diagonal entries are independent integral coordinates. The quotient
diagonal lattice is \(\Z^r\) in odd dimension. In even dimension it is
\(\Z^r\cup(\Z+\tfrac12)^r\), since all \(\gamma\pm u_i\) must be integers.
Its lattice density \(\rho_{\rm diag}\) is therefore \(1\) in odd dimension and \(2\)
in even dimension. At fixed \(u\), \(\gamma\) has spacing one.

Write \(v=\max_i|u_i|\). The scalar interval over \(u\) in the entrywise unit box has
length \(2(1-v)\) for \(v\le1\), and the interval over \(tu\) has length \(2(1-tv)\)
for \(0\le t\le1\). In \eqref{eq:slice} the off-diagonal coordinates give \(2^{m-r}\)
and the cube shell in \(u\) has density \(\rho_{\rm diag}\,r2^r v^{r-1}\,dv\).
Dividing by \(2^{m+2}\) gives
\[
 f_B(t)=\rho_{\rm diag}\,r\int_0^1v^{r-1}(1-v)(1-tv)\,dv
 =\frac{\rho_{\rm diag}}{r+1}\left(1-\frac r{r+2}t\right).
\]
Formula \eqref{eq:affine} gives the first three entrywise coefficients.

\head{Split dimension four}
The determinant form on \(M_2\) has polar matrix \(J_4\) up to a signed
permutation of the four basis coordinates. Conjugating a matrix by the
corresponding signed permutation matrix permutes its entries and changes some
signs. It therefore preserves both the integral matrix lattice and the entrywise
box. In the determinant model the integral matrices are
\[
 tI_4+I_2\otimes\begin{pmatrix}0&e\\ f&p\end{pmatrix}
 +\begin{pmatrix}0&g\\ h&q\end{pmatrix}\otimes I_2,
 \qquad t,p,q,e,f,g,h\in\Z.
\]
Indeed left and right multiplication on \(M_2\) give the two tensor summands and
span the six-dimensional quotient by scalars. The entries recover \(t\), the four
off-diagonal parameters, and the diagonal differences \(p,q\). Thus the quotient
lattice is the product of the triples \((p,e,f)\) and \((q,g,h)\) in its two
factors of type \(A_1\). Their coordinate boxes give \(\bigl(4\zeta(2)/\zeta(3)\bigr)^2\) by \eqref{eq:product}.

For the entrywise box, let \(\ell=(a,b)\) and \(\ell'=(a',b')\) range independently
over primitive integer slopes modulo common sign. Write \(j=\max(|a|,|b|)\),
\(k=\max(|a'|,|b'|)\), \(x=\min(|a|,|b|)/j\) and \(y=\min(|a'|,|b'|)/k\).
Then
\begin{equation}\label{eq:gso4}
 \begin{aligned}
 K_4&=\sum_{\ell,\ell'}\frac{w_4(\ell,\ell')}{j^3k^3},\\
 w_4(\ell,\ell')&=\begin{cases}
 \displaystyle\frac23-\frac{x+y}{6},&\text{maxima on the same side},\\[5pt]
 \displaystyle\frac{(1-x)(3-x)+(1-y)(3-y)}{6(1-xy)},&\text{otherwise}.
 \end{cases}
 \end{aligned}
\end{equation}
The maxima are on the same side if both occur in the first coordinates or both in
the second. At a tie the first expression applies. Indeed the four diagonal offsets
from \(t\) are \(0,p,q,p+q\), of range \(|p|+|q|\). Integrating the two scalar
intervals on a slope section and rescaling the absolute diagonal parameters by
\(j/(2H),k/(2H)\) gives
\[
 w_4(\ell,\ell')=4\int_0^1\!\int_0^1
 (1-|a|u/j-|a'|v/k)_+(1-|b|u/j-|b'|v/k)_+\,du\,dv ,
\]
where \(w_+=\max(w,0)\). Each triple contributes its slope-height factor
\(j^{-3}\) or \(k^{-3}\). For maxima on the same side the support is the triangle
\(u+v\le1\). In the other case its vertices are \((0,0),(1,0)\),
\(\big((1-y)/(1-xy),(1-x)/(1-xy)\big),(0,1)\) after interchanging the two
matrices if necessary. Integrating over the triangles from the origin gives
\eqref{eq:gso4}. The weights are bounded and \(\sum_\ell j^{-3}=4\zeta(2)/\zeta(3)\), so
the series converges absolutely. The error is \(O(H^7\log(2H))\) by \eqref{eq:error}.
Finally, \(\gso_1\) and split \(\gso_2\) consist of scalar and diagonal matrices,
respectively. Their entrywise counts are exactly \(L^2\) and \(L^4\).

\subsubsection{The adjoint representation}
For the exceptional types we use the adjoint representation on the split algebra
with Chevalley basis \((h_1,\ldots,h_r,(e_\alpha)_{\alpha\in\Phi})\), where
\(h_i\) are the simple coroots. The entry bound is on the \(n\times n\) matrix of
\(\ad x\) in this basis. The volume calculation separates into one Cartan polytope
and one interval for each root. To record the widths of these intervals we use three
kinds of entries of \(\ad e_\alpha\). For \(\alpha=\sum c_i\alpha_i\) write
\(h_\alpha=\sum c_i^\vee h_i\). For \(\alpha+\beta\in\Phi\) put
\(p_{\alpha,\beta}=\max\{j\ge0:\beta-j\alpha\in\Phi\}\). Set
\begin{equation}\label{eq:rootmax}
 b_\alpha=\max\left\{\max_i|\alpha(h_i)|,\ \max_i|c_i^\vee|,\
 \max_{\beta:\alpha+\beta\in\Phi}(1+p_{\alpha,\beta})\right\},
\end{equation}
with the last maximum omitted when its index set is empty. Let \(\sum m_i\alpha_i\)
be the highest root, \(W_\Phi\) the Weyl group, and
\[
 v_\Phi=\frac{|W_\Phi|}{r!\prod_i m_i}.
\]

\begin{proposition}\label{pr:adjoint}
For the full integral adjoint matrix lattice just specified,
\begin{equation}\label{eq:adjoint}
 K_{\rm ad}=\frac{4v_\Phi}{2^r\prod_{\alpha\in\Phi}b_\alpha}\,\frac{\zeta(n-1)}{\zeta(n)}.
\end{equation}
\end{proposition}
\begin{proof}
Write \(x=h+\sum_\alpha t_\alpha e_\alpha\). Relative to the Cartan weight
decomposition, the matrices \(\ad e_\alpha\) for distinct roots have disjoint
entry supports. Indeed the row weight minus the column weight of any nonzero
entry equals \(\alpha\), also when one weight is zero. The Cartan part is
diagonal with entries \(\beta(h)\) on the root lines and zero on the Cartan.

In \(\ad e_\alpha\) the columns belonging to \(h_i\) have entry \(-\alpha(h_i)\)
and the column belonging to \(e_{-\alpha}\) has entries \(\pm c_i^\vee\).
The remaining nonzero entries have absolute values \(1+p_{\alpha,\beta}\) by
the Chevalley basis relations \cite[\S1]{Steinberg}. Thus the root coordinate has
interval \(|t_\alpha|\le b_\alpha^{-1}\). The coroot \(h_\alpha\) is primitive
in the coroot lattice, being Weyl-conjugate to a simple coroot. Its coefficients
have greatest common divisor one, so integrality of the matrix forces
\(t_\alpha\in\Z\) and this condition suffices for its root block.

The Cartan lattice of integral adjoint matrices is
\(\{h:\alpha_i(h)\in\Z,\ 1\le i\le r\}\). Thus the coordinates \(u_i=\alpha_i(h)\)
have lattice covolume one. In the dominant chamber the inequalities
\(|\beta(h)|\le1\) for all roots reduce to \(u_i\ge0\) and \(\sum m_i u_i\le1\).
Its simplex volume is \((r!\prod_i m_i)^{-1}\). Since the Weyl group preserves the
lattice and the chambers partition the body up to measure zero, its Cartan volume
is \(v_\Phi\). Multiplying by \(\prod_\alpha(2/b_\alpha)\) and using
\(n=r+|\Phi|\) in \eqref{eq:representation} proves the formula. In particular the
lattice here includes the entire coweight lattice in its Cartan part.
\end{proof}

\subsubsection{Evaluation in the five exceptional types}
Let \(a_j\) count the roots of both signs with \(b_\alpha=j\), so that
\(\prod_{\alpha\in\Phi}b_\alpha=\prod_jj^{a_j}\). The table gives the data
in \eqref{eq:adjoint} and the resulting coefficients of \(L^n\).
\begin{center}
\small\setlength{\tabcolsep}{5pt}\renewcommand{\arraystretch}{1.65}
\begin{tabular}{@{}lrrrrrrrrl@{}}
\toprule
Type& \(n\)&\(v_\Phi\)&\(a_1\)&\(a_2\)&\(a_3\)&\(a_4\)&\(a_5\)&\(a_6\)&\(K_{\rm ad}\)\\
\midrule
\(G_2\)&14&1&0&2&10&0&0&0&\(\dfrac{\zeta(13)}{2^2\,3^{10}\,\zeta(14)}\)\\\addlinespace[4pt]
\(F_4\)&52&1&0&38&4&6&0&0&\(\dfrac{\zeta(51)}{2^{52}\,3^4\,\zeta(52)}\)\\\addlinespace[4pt]
\(E_6\)&78&3&38&30&4&0&0&0&\(\dfrac{\zeta(77)}{2^{34}\,3^3\,\zeta(78)}\)\\\addlinespace[4pt]
\(E_7\)&133&2&54&50&16&6&0&0&\(\dfrac{\zeta(132)}{2^{66}\,3^{16}\,\zeta(133)}\)\\\addlinespace[4pt]
\(E_8\)&248&1&72&76&40&30&12&10&\(\dfrac{\zeta(247)}{2^{152}\,3^{50}\,5^{12}\,\zeta(248)}\)\\
\bottomrule
\end{tabular}
\end{center}
For example, the \(G_2\) lattice has the asymptotic
\[
 N_{\Lambda_{\rm ad},B_{\rm ad}}(H)
 =\frac{\zeta(13)}{236196\,\zeta(14)}L^{14}+O(H^{14-1/3}).
\]
With the short simple root first, its positive roots in simple-root coordinates
are \((1,0),(0,1),(1,1),(2,1),(3,1),(3,2)\). Formula \eqref{eq:rootmax}
gives successive maxima \(3,3,3,3,3,2\). The highest root has marks \((3,2)\)
and \(|W_\Phi|=12\), so \(v_\Phi=1\).

For completeness we give the finite enumeration producing the other rows.
It requires only the Cartan matrix and also specifies the root numbering. Let
\(A_{ij}=\alpha_j(h_i)\), with diagonal entries two. For \(F_4\) use a chain
with squared lengths proportional to \((2,2,1,1)\). For \(E_r\) take the simply
laced chain \(1,\ldots,r-1\) and attach \(r\) to vertex \(3\). For \(G_2\) the
length factors are \((1,3)\). Write \(d_i\) for these factors, taking \(d_i=1\)
in the simply laced cases. On each edge \(i,j\) put
\(A_{ij}=-\max(1,d_j/d_i)\). All nonadjacent off-diagonal entries are zero.
The following recurrence reproduces every multiplicity in the table. The accompanying
Python supplement \texttt{exceptional\_entry\_boxes.py} is also attached to this PDF.
Running it prints the multiplicities, highest root and rational factor \(K_{\rm ad}\,\zeta(n)/\zeta(n-1)\)
for each of the five types.
\begin{enumerate}
\item Start with the unit columns \(c=e_i\) and close the set under
\(c\mapsto c-(Ac)_i e_i\) for all \(i\). This gives precisely the root columns,
since every root is Weyl-conjugate to a simple root.
\item For each column put \(q(c)=\tfrac12\sum_i d_i c_i(Ac)_i\).
The coroot has coefficients \(d_i c_i/q(c)\), by the formula
\(\alpha^\vee=2\alpha/(\alpha,\alpha)\). Initialize \(b_\alpha\) to the larger
of \(\norm{Ac}\) and \(\max_i|d_i c_i/q(c)|\).
\item For each root column \(b\) with \(b+c\) a root, let \(j\) be the first
positive integer for which \(b-jc\) is absent. Replace \(b_\alpha\) by
\(\max(b_\alpha,j)\). Increment \(a_{b_\alpha}\) once when all such \(b\) have
been considered.
\end{enumerate}
In the chosen numbering the highest root columns for \(F_4,E_6,E_7,E_8\) are
respectively \((2,3,4,2)\), \((1,2,3,2,1,2)\), \((2,3,4,3,2,1,2)\) and
\((2,4,6,5,4,3,2,3)\). Their Weyl group orders are \(1152,51840,2903040,696729600\)
\cite[Ch.~VI, plates]{Bourbaki}. These data give the displayed \(v_\Phi\). As a check,
in the simply laced cases all root-string entries in \eqref{eq:rootmax} have absolute
value one. For a positive nonsimple root the maximum is therefore its largest
simple-root coefficient. For a simple root it is two.

If the lattice is instead the integer span of the chosen Chevalley basis, the Cartan
lattice is the coroot lattice. Its index in the coweight lattice is \(|\det A|\),
and the root coordinates have the same lattice in both models. For the same adjoint
entry box the Chevalley lattice coefficient is \(K_{\rm ad}/|\det A|\).
The affine Weyl alcove identity \(|W_\Phi|=r!(\prod_i m_i)|P_\Phi/Q_\Phi|\),
where \(P_\Phi,Q_\Phi\) are the weight and root lattices, gives
\(v_\Phi=|P_\Phi/Q_\Phi|=|\det A|\) \cite[Ch.~VI]{Bourbaki}.
Consequently the Chevalley lattice gives, for the same adjoint entry box,
\[
 K_{\rm Ch}=\frac{4}{2^r\prod_{\alpha\in\Phi}b_\alpha}\,\frac{\zeta(n-1)}{\zeta(n)}.
\]
In the order \(G_2,F_4,E_6,E_7,E_8\) the indices are \(1,1,3,2,1\).

\subsection*{AI disclosure} The authors thank ChatGPT for assistance with proof-strategy discussions, exploratory computations, and the editing of this manuscript.

\end{document}